\documentclass[preprint,10pt]{elsarticle}
\usepackage{amssymb}
\usepackage{graphicx}
\usepackage{amsmath}
\usepackage{amsthm}
\usepackage{lineno}
\biboptions{comma,square}
\biboptions{}
\numberwithin{equation}{section}
\usepackage{a4wide}

\def\qed   {\hfill $\scriptstyle \square$}

\def\RR    {{R\!\!\!\!\!I~}}

\newcommand \Scal {\mathcal S}
\newcommand \Rcal {\mathcal R}
\newcommand \CC {\mathcal C} 
\def\ld    {\lambda}

\def\ZZ    {\bold Z}

\def\LL    {{\mathcal L}}
\def\NN    {{\mathbf N}}

\newcommand \WW  {\mathcal W}

\def\pt    {\partial}

\def\RR    {{R\!\!\!\!\!I~}}

\def\pt   {\partial}

\def\NN    {{\mathbf N}}

\newcommand{\rf}[1]{(\ref{#1})}
\makeatletter
\def\ps@pprintTitle{%
     \let\@oddhead\@empty
     \let\@evenhead\@empty
     \def\@oddfoot{\footnotesize\itshape
      \ifx\@journal\@empty 
       \else\@journal\fi\hfill {\bf Journal of Computational Physics (2011)}}%
     \let\@evenfoot\@oddfoot}
\makeatother

\begin{document}

\begin{frontmatter}

\title
{A Godunov-type method for the shallow water equations
with discontinuous topography 
 in the resonant regime} 

\author%[P.G. L{\tiny e}Floch and M.D. Thanh]
{  Philippe G.  LeFloch$^1$ and  Mai Duc Thanh$^2$}

\address{
$^1$
%%%%%%%%%% SPELLING:             LeFLOCH     or   LeFloch 
%%%%
Laboratoire Jacques-Louis Lions \& Centre National de la Recherche
Scientifique
\\
Universit\'e Pierre et Marie Curie (Paris 6), 
4 Place Jussieu, 75252 Paris,
France. \\
Email: pgLeFloch@gmail.com. Web: philippelefloch.wordpress.com.}

\address{
$^2$
%Mai Duc Thanh \newline
Department of Mathematics, International University, Quarter 6, 
Linh Trung Ward
\\
Thu Duc District, Ho Chi Minh City, Vietnam. 
Email: mdThanh@hcmiu.edu.vn}

%\subjclass[2000]{35L65, 76N10, 76L05}

\begin{abstract} We investigate the Riemann problem for the shallow water equations
with variable and (possibly) discontinuous topography and provide a complete description of 
the properties of its solutions:  
existence; uniqueness in the non-resonant regime; multiple solutions in the resonant regime.
This analysis leads us to a numerical algorithm that provides one with a Riemann solver.
Next, we introduce a Godunov-type scheme based on this Riemann solver, which is
well-balanced and of quasi-conservative form. Finally, we present numerical 
experiments which demonstrate the convergence of the proposed scheme even in the resonance regime,
except in the limiting situation when Riemann data precisely belong to the resonance hypersurface.
\end{abstract}


\end{frontmatter}

\newtheorem{theorem}{Theorem}[section]
\newtheorem{lemma}[theorem]{Lemma}
\newtheorem{proposition}[theorem]{Proposition}
\newtheorem{definition}[theorem]{Definition}
\newtheorem{corollary}[theorem]{Corollary}
\newdefinition{rmk}{Remark}

\section{Introduction}

In this paper we design a Godunov-type scheme for the numerical approximation
 of weak solutions
to the initial-value problem associated with the
shallow water equations with variable topography, i.e. 
\begin{equation}
\begin{array}{ll}
&\pt_th + \pt_x(hu)  =  0,
\\
&\pt_t(hu) + \pt_x \Big( h(u^2 + {g h \over 2}) \Big)  + gh \, \pt_x a = 0,
\\
\end{array} \label{1.1}
\end{equation}
where the height of
the water from the bottom to the surface, denoted by $h$, and the fluid velocity $u$ are the main unknowns. 
Here, $g$ is the so-called gravity constant, and $a=a(x)$ (with $x\in\RR$) is the height of the
bottom from a given level. 

In \cite{LeFloch89}, LeFloch pointed out that by supplementing balance laws like \eqref{1.1} 
with the additional equation
\begin{equation}
\pt_t a=0,
\label{1.2}
\end{equation}
the set of equations \rf{1.1}--\rf{1.2} can be regarded as a {\sl nonlinear hyperbolic system
in nonconservative form} 
and tackled with the theory introduced by Dal~Maso, LeFloch, and Murat \cite{DLM}
and developed by LeFloch and collaborators \cite{LeFloch88,LeFlochLiu,HL,CPMP}.
As is well-known, the system \rf{1.1}--\rf{1.2} is hyperbolic but {\sl not strictly hyperbolic} since
characteristic speeds may coincide on certain hypersurfaces.
Non-strictly hyperbolic systems have been extensively studied in the literature.
See \cite{LeFloch89, MarchesinPaes-Leme, LeFlochThanh03.2, Thanh09} for the model of fluid flows in a nozzle with variable cross-section, \cite{IsaacsonTemple95, IsaacsonTemple92, GoatinLeFloch, AndrianovWarnecke, AndrianovWarnecke2} for other models.

On the other hand, 
the discretization of source terms in nonlinear hyperbolic balance laws like the shallow water equations
was pioneered by Greenberg and Leroux \cite{GreenbergLeroux}. 
We built here on this paper as well as the follow-up work 
\cite{AlcrudoBenkhaldoun,Bouchutbook04,Bernettietal,ChinnayyaLeRouxSeguin,Gosse00,GreenbergLerouxBarailleNoussair,JinWen1,JinWen2}.
Since the system \eqref{1.1} is also related to the class of two-phase models, we are also motivated by the existing
research work on the discretization of two-phase flow models 
\cite{CoqueletalM2AN2009,AndrianovWarnecke2,GallouetHerardSeguin}. 
In particular, recall that well-balanced schemes for shallow water equations were constructed first in 
\cite{GallouetHerardSeguin,GreenbergLeroux,JinWen1, JinWen2}.
In addition, the discretization of nonconservative hyperbolic systems
and of systems with source terms  attracted a lot of attention in recent years.
We refer to \cite{AudusseBouchutBristeauKleinPerthame,BotchorishviliPerthameVasseur, BotchorishviliPironneau, Gosse00,GreenbergLerouxBarailleNoussair} for a single conservation law with source term
and to
\cite{ JinWen1, JinWen2,KroenerThanh05, KroenerLeFlochThanh} for fluid flows in a nozzle with variable cross-section. Well-balanced schemes for multi-phase flows and other models were studied in
\cite{CoqueletalM2AN2009,Bouchutbook04, SaurelAbgrall99JCP, ThanhKroenerNam}.

The Godunov scheme is based on an exact or approximate Riemann solver and, consequently, it is necessary that 
sufficient information be available on the existence and properties of all solutions to \rf{1.1}--\rf{1.3}.
Recall that the Riemann problem is a Cauchy problem
with piecewise constant initial data of the form
\begin{equation}
(h,u,a)(x,0)= \left\{\begin{array}{ll}(h_L,u_L,a_L),\quad &x< 0,\\
(h_R,u_R,a_R),\quad &x>0. \end{array}\right. \label{1.3}
\end{equation}
One main task in the present paper will be to revisit the construction of
solutions to the Riemann problem for \rf{1.1}--\rf{1.3} in order to arrive at a definite algorithm for their numerical computation.
We will show that a unique solution exists within a large domain of initial data,
and will precisely identify the domains where multiple solutions are available, 
by providing necessary and sufficient conditions for the existence of multiple (up to three) solutions.

Recall that, in LeFloch and Thanh \cite{LeFlochThanh07}, a first investigation of
the Riemann problem for the shallow water equations  
was performed. The present paper provides a very significant improvement in that
a complete description of a Riemann solver is now obtained for the first time.
We are able to distinguish between cases of existence, uniqueness, and multiplicity of solutions. 

We refer the reader to \cite{AlcrudoBenkhaldoun,Bernettietal, RosattiBegnudelli} for partial or alternative
approaches to the Riemann problem.
On the other hand, Chinnayya, LeRoux, and Seguin  \cite{ChinnayyaLeRouxSeguin} introduced 
a Godunov method for \rf{1.1} based 
on a Riemann solver determined by ``continuation": they start 
their construction by assuming that the bottom topography is flat and then extend it to a non-flat topography.
Their method allows them to construct solutions within the regime
where the nonlinear characteristic fields are separated by the linearly degenerate one, that is,
the case where one wave speed is negative and the other positive.
 In this regime, the Riemann solution with non-flat bottom can be obtained from the
the wave curves associated with the fast wave family, only, and these two waves are
separated by a stationary wave. 
In particular, the total number of waves in Riemann solutions is exactly
the number of characteristic families.
On the other hand,
when more general Riemann data are considered and lie in the other two strictly hyperbolic regions,
say of ``cross-region" type, we find it hard to apply this method.
There is always a major difference between the ``flat-bottom" and ``non-flat-bottom" cases, since
in one case the system is strictly hyperbolic, the other case the system is not strictly hyperbolic.
Indeed, in the non-flat-bottom case, new wave curves arise that
replace more standard wave curves.
The total number of waves in a solution can possibly be {\sl larger}
than the number of characteristic fields as waves associated with a given family can be repeated.
The appearance of such new wave curves cannot be obtained by a continuation argument.

Our objective in the present work
is to present a numerical algorithm that provides an explicit construction of a Riemann solver for the shallow water equations in the resonant or non-resonant regimes.
This solver then provides us with solutions to
local Riemann problems which we can use
to design a Godunov-type scheme.
We also provide extensive numerical experiments and, within strictly hyperbolic regions 
of the phase space, our 
tests demonstrate that the proposed scheme converge to the expected solution. In resonance regions, 
we also observe convergence, except when the Godunov scheme takes some
 values on the resonance hypersurfaces.
The Godunov scheme proposed in the present paper for \rf{1.1}--\rf{1.2} turns out to be well-balanced and captures exactly stationary waves.

This paper is organized as follows. In Section~2, we recall basic facts about
the system \rf{1.1}--\rf{1.2} and provide the computing algorithm for stationary contact waves.
In Section~3, we revisit the construction of solutions to the Riemann problem, and
present a new approach based on a ``gluing" technique involving different solution structures.
With this technique, we establish the existence of solutions for arbitrary large initial data.
In Section~4 we then present our computing strategy leading to a Riemann solver and,
finally, are in a position to design a corresponding Godunov scheme. Section~5 is devoted to 
numerical experiments when data belong to strictly hyperbolic domains; in particular, 
we estimate the  numerical errors and observe convergence when the mesh size tends to zero.
Finally, Section~6 is devoted to numerical tests here data belong to resonance regions, and precisely identify
regimes of convergence or non-convergence.

%*****************************************************************************************

\section{Shallow water equations}

\subsection{Wave curves}

Introducing the dependent variable
$(h,u,a)=(h,u,a)(x,t)$, the Jacobian matrix of the system  \rf{1.1}--\rf{1.2} admits three real eigenvalues, i.e.
\begin{equation}
\ld_1(U):=u-\sqrt{gh}<\ld_2(U):=u+\sqrt{gh}, \quad \ld_3(U):=0,
\label{2.1}
\end{equation}
together with the corresponding  eigenvectors:
\begin{equation}
r_1(U):=\left(\begin{matrix}h\\
-\sqrt{gh}\\
0\end{matrix}\right),\quad
r_2(U):=\left(\begin{matrix}h\\
\sqrt{gh}\\
0\end{matrix}\right),\quad r_3(U):=\left(\begin{matrix}gh\\
-gu\\
u^2-gh\end{matrix}\right),
\label{2.2}
\end{equation}
so that the system \rf{1.1}--\rf{1.2} is hyperbolic, but not strictly hyperbolic.
More precisely,  the first and the third characteristic speeds coincide,
$$
(\ld_1(U),r_1(U))=(l_3(U),r_3(U)),
$$
 on the hypersurface
\begin{equation}
\CC_+:= \Big\{ (h,u,a) | \quad u=\sqrt{gh} \Big\}.
\label{2.3}
\end{equation}
The second and the third characteristic fields coincide,
$$
(\ld_2(U),r_2(U))=(\ld_3(U),r_3(U)),
$$
on the hypersurface
\begin{equation}
\CC_-:= \Big\{(h,u,a) | \quad u=-\sqrt{gh} \Big\}.
\label{2.4}
\end{equation}

\begin{center}
\begin{figure}[h]
  \includegraphics[width=8cm]{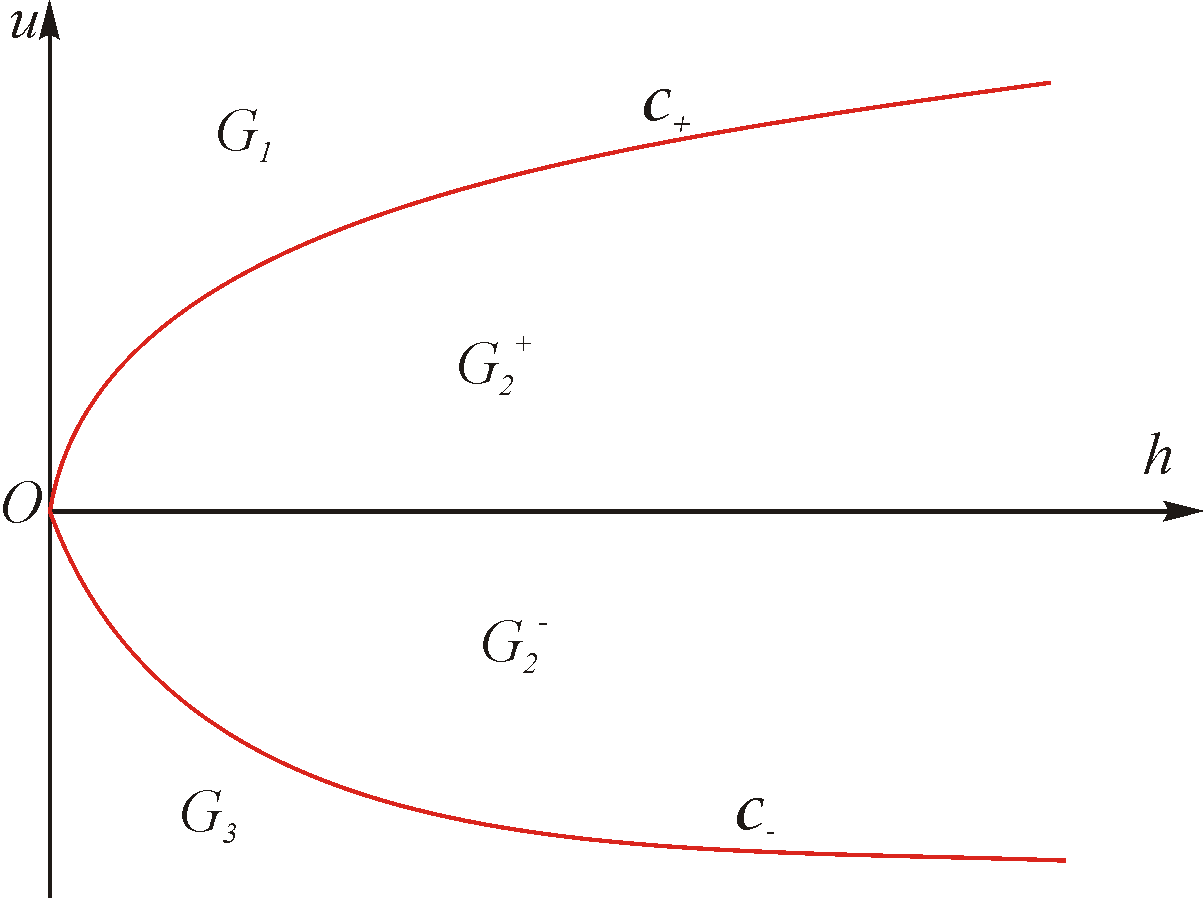}\\
  \caption{Phase domain in the $(h,u)$-plane.}\label{phaseplane}
\end{figure}
\end{center}

In the following, we introduce $\CC:=\CC_+\cup \CC_-$,
and we refer to the sets $\CC_\pm$  as the {\it resonance hypersurfaces},
which
separate the phase space in the $(h,u,a)$-variable
into three sub-domains
\begin{equation}
\begin{array}{ll}&G_1:=
\Big\{(h,u,a) \in \RR_+\times\RR\times\RR_+ |\quad \ld_1(U)>\ld_3(U) \Big\},
\\
&G_2:= \Big\{(h,u,a) \in \RR_+\times\RR\times\RR_+ | \quad \ld_2(U)>\ld_3(U)>\ld_1(U) \Big\},
\\
&G_3:= \Big\{(h,u,a)\in \RR_+\times\RR\times\RR_+  |\quad \ld_3(U)>\ld_2(U) \Big\},
\end{array}\label{2.5}
\end{equation}
in which the system is strictly hyperbolic.
It is convenient to also set
$$
G_2^+:=\{(h,u,a) \in G_2, u\ge 0\},\qquad G_2^-:=\{(h,u,a) \in G_2, u\le 0\}.
$$
Observe that for $u\ge 0$, the set $G_1$ is referred to as the domain of supercritical flows, where the Froude number
$$
Fr:={u\over\sqrt{gh}}
$$
is larger than 1, and the set $G_2^+$ as the domain of subcritical flows, where $Fr<1$. This concept may similarly be extended to the case $u<0$.
See Figure \ref{phaseplane}.

One easily checked that the first and second
characteristic fields $(\ld_1,r_1)$, $(\ld_2,r_2)$ are genuinely
nonlinear, while  the  third characteristic field $(\ld_3,r_3)$ is linearly degenerate.

%----------------------------------------------------------------------------------------

As discussed in \cite{LeFlochThanh07}, across a discontinuity there are two possibilities:
\begin{itemize}
\item[(i)] either the bottom height $a$ remains constant,
\item[(ii)] or the discontinuity is stationary (i.e.~propagates with zero speed).
\end{itemize}
In the first the case (i), the system \rf{1.1}--\rf{1.2} reduces to the standard
shallow water equations with flat bottom.
We can determine the $i$-shock curve $\Scal_i(U_0)$ starting
from  a left-hand state $U_0$ and consisting of all right-hand states
$U$ that can be connected to $U_0$ by a Lax shock associated with
the first characteristic field ($ i=1,2$):
\begin{equation}
\Scal_i(U_0) :\quad \Psi_i(U;U_0):=u-u_0\pm  \sqrt{g\over
2}(h-h_0)\sqrt{\Big({1\over h}+{1\over h_0}\Big)}=0.
\label{2.6}
\end{equation}
where $U=(h,u)$, $h> h_0$ for $i=1$ and $\quad h<h_0$ for $i=2$.
We also define the backward $i$-shock curve $\Scal_i^B(U_0)$
starting from a right-hand state $U_0$ and consisting of all left-hand
states $U$ that can be connected to $U_0$ by a Lax shock
associated with the first characteristic field $(i=1,2$):
\begin{equation}
\Scal_i^B(U_0) :\quad \Phi_i(U;U_0):=u-u_0\pm  \sqrt{g\over
2}(h-h_0)\sqrt{\Big({1\over h}+{1\over h_0}\Big)}=0,
\label{2.7}
\end{equation}
where $U=(h,u)$, $h< h_0$ for $i=1$ and $\quad h>h_0$ for $i=2$.

It is interesting that the shock speed of the nonlinear characteristic fields may coincide with the speed of stationary contact waves. The following lemma is easily checked. 

\begin{lemma} \label{lem31}   Consider the projection on the $(h,u)$-plan. To every $U=(h,u)\in G_1$ there exists exactly one point $ U^{\#}\in\Scal_1(U)\cap G_2^+$ such that the $1$-shock speed $\bar\ld_1(U,U^{\#})=0$.
The state $U^{\#} =(h^{\#},u^{\#})$ is defined by
$$
h^{\#}=\frac{-h+\sqrt{h^2+8hu^2/g}}{2},\quad u^{\#}=\frac{uh}{h^{\#}}.
$$
Moreover, for any $V\in \Scal_1(U)$, the shock speed $\bar\ld_1(U,V)>0$ if and only if $V$ is located above $U^{\#}$ on $\Scal_1(U)$.
\end{lemma}

It is also well-known that the bottom height $a$ remains constant through rarefaction fans.
 The forward rarefaction curve $\Rcal_i(U_0)$
starting from a given left-hand state $U_0$ and consisting of all the
right-hand states $U$ that can be connected to $U_0$ by a
rarefaction wave associate with the first characteristic field as
\begin{equation}
\Rcal_i(U_0):\quad
 \Psi_i(U;U_0)=u-u_0\pm 2\sqrt{g}(\sqrt{h}-\sqrt{h_0})=0,\quad i=1,2,\label{2.8}
\end{equation}
where $U=(h,u), h\le  h_0$ for $i=1$ and $h\ge h_0$ for $i=2$.
Given a right-hand state $U_0$, the backward $i$-rarefaction
curve $\Rcal_i^B(U_0)$ consisting of all left-hand states $U$ that
can be connected to $U_0$ by a rarefaction wave associated with the
first characteristic field reads $(i=1,2$)
\begin{equation}
\Rcal_i^B(U_0):\quad
 \Phi_i(U;U_0)=u-u_0\pm 2\sqrt{g}(\sqrt{h}-\sqrt{h_0})=0,
\label{2.9}
\end{equation}
where $U=(h,u), h\ge  h_0$ for $i=1$ and $h\le h_0$ for $i=2$.

Finally, we define the forward and backward wave curves in the $(h,u)$-plane
$(i=1,2$):
\begin{equation}
\begin{array}{ll}&\WW_i(U_0):=\Scal_i(U_0)\cup\Rcal_i(U_0) =\{U \ | \ \Psi_i(U;U_0)=0\},\\
&\WW_i^B(U_0):=\Scal_i^B(U_0)\cup\Rcal_i^B(U_0) =\{U \ | \ \Phi_i(U;U_0)=0\}.
\end{array} \label{2.10}
\end{equation}
It is checked in \cite{LeFlochThanh07} that the wave curves $\WW_1(U_0)$
and $\WW_1^B(U_0) $
 parameterized as $h\mapsto u=u(h), h>0,$ are strictly convex and strictly decreasing
functions. The wave curve $\WW_2(U_0)$ and
$\WW_2^B(U_0)$  being parameterized as
$h\mapsto u=u(h), h>0,$ are strictly concave and strictly
decreasing functions.

In the case (ii),
the discontinuity satisfies the jump relations
\begin{equation}
\begin{array}{ll}&[hu]  =  0,
\\
&[{u^2\over 2} + g(h+a)]  = 0,
\\
\end{array}\label{2.11}
\end{equation}
which determine the stationary-wave curve (parameterized with $h$):
\begin{equation}
\begin{array}{ll}\WW_3(U_0):\quad &u=u(h)={h_0u_0\over h},\\
& a=a(h)=a_0+{u_0^2-u^2\over 2g}+h_0-h. \end{array} \label{2.12}
\end{equation}
The projection of the wave curve $\WW_3(U_0)$  in the $(h,u)$-plane can be parameterized as $h\mapsto
u=u(h), h>0,$ which is a strictly convex and strictly decreasing
function for $u_0>0$ and strictly concave and strictly increasing
function for $u_0<0$.

The above arguments show that {\it the $a$-component of Riemann solutions may
change only across a stationary wave.}
This property will be
 important later when designing the discretization of the source terms.

\subsection{Properties of stationary contacts}

\begin{center}
\begin{figure}
  \includegraphics[width=9cm]{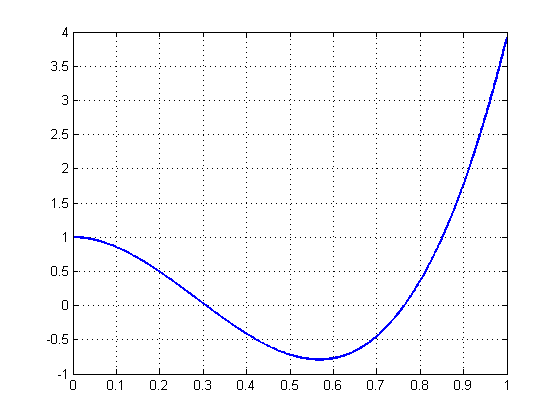}\\
  \caption{Graph of the function $\varphi=\varphi(h), h\ge 0$ defined by \rf{2.13}
with $g=9.8, a_0=1,h_0=1, u_0=1$ and $a=1.2$.
The function $\varphi$ admits two zeros in the interval $(0,1)$.}\label{varphi}
\end{figure}
\end{center}

Given a state $U_0=(h_0,u_0,a_0)$ and another bottom level $a\ne a_0$, we let
 $U=(h,u,a)$ be the corresponding right-hand state of the stationary contact issuing
from the given left-hand state $U_0$. We now determine $h, u$ in terms of $U_0,a$,
as follows. Substituting $u=h_0u_0/h$ from the first equation of \rf{2.12} to the second equation of \rf{2.12}, we obtain
$$
a_0+{1\over 2g}\left(u_0^2- \left({h_0u_0\over h}\right)^2\right)+h_0-h=a.
$$
Multiplying both sides of the last equation by $2gh^2$, and then re-arranging terms,
we find that $h>0$ is a root of the nonlinear equation
\begin{equation}
\aligned
\varphi(h)
& = \varphi(U_0,a;h)
\\
& :=2gh^3 + (2g(a-a_0-h_0)-u_0^2)h^2 + h_0^2u_0^2=0.
\endaligned
 \label{2.13}
\end{equation}
We easily check
$$
\begin{array}{ll}&\varphi(0)=h_0^2u_0^2\ge 0,\\
&\varphi'(h)=6gh^2+2(2g(a-a_0-h_0)-u_0^2)h,\\
&\varphi''(h)=12gh+2(2g(a-a_0-h_0)-u_0^2),
 \end{array}$$
so that
\begin{equation}
\aligned
\varphi'(h) = 0 \quad&\textrm{iff}\quad h=0
\\
&  \textrm{or}\quad  h=h_*=h_*(U_0,a):=\frac{u_0^2+2g(a_0+h_0-a)}{3g}.
\endaligned
  \label{2.14}
\end{equation}
 If $h_*(U_0,a)<0$, or $ a>a_0+h_0 +{u_0^2\over 2g},$
then $\varphi'(h)>0$ for $h>0$. Since $\varphi(0)=h_0^2u_0^2\ge 0$, there is no root for \rf{2.13} if \rf{2.14} holds. Otherwise, if
$$
  a\le a_0+h_0 +{u_0^2\over 2g},
$$
then $\varphi'>0$ for $h>h_*$ and $\varphi'(h)<0$ for $0<h<h_*$. In this case, $\varphi$ admits a zero $h>0$, and in this case it has two zeros, iff
$$
\varphi_{\min}:=\varphi(h_*)=-gh_*^3 + h_0^2u_0^2\le 0,
$$
 or
\begin{equation}
h_*(U_0,a)\ge h_{\min}(U_0) :=\Big({h_0^2u_0^2\over g}\Big)^{1/3},
  \label{2.15}
\end{equation}
where $h_*$ is defined by \rf{2.14}. It is easy to check that \rf{2.15} holds if and only if
\begin{equation}
\begin{array}{ll}
  a\le a_{\max}(U_0)&:=a_0+h_0 +{u_0^2\over 2g}-{3\over 2g^{1/3}}(h_0u_0)^{2/3}\\
 &=  a_0 + {1\over 2g}\big((gh_0)^{1/3}-u_0^{2/3}\big)^2(2(gh_0)^{1/3}+u_0^{2/3}).
 \end{array}
 \label{2.16}
\end{equation}
Observe that \rf{2.16} implies $a_{\max}(U_0)\ge a_0$ and the equality holds only if $(h_0,u_0)$ belongs to the surfaces $\CC_\pm$.
Whenever \rf{2.16} is fulfilled, the function $\varphi$ in \rf{2.13} admits two roots denoted by  $h_1(a)\le h_2(a)$ satisfying $h_1(a)\le h_*\le h_2(a)$. Moreover, if the inequality in \rf{2.16} is strict, i.e., $a<a_{\max}(U_0)$, then these two roots are distinct: $h_1(a)< h_* < h_2(a)$. See Figure \ref{varphi}.

Thus, we arrive at the following lemma.

\begin{lemma} \label{lem21} Given $U_0=(h_0,u_0,a_0)$ and a bottom level $a\ne a_0$. The following conclusions holds.
\begin{itemize}
\item[(i)]
$
a_{\max}(U_0)\ge a_0,\quad a_{\max}(U_0)=a_0\quad\textrm{ if and only if}\quad (h_0,u_0)\in\CC.
$
\item[(ii)] The nonlinear equation \rf{2.13} admits a root if and only if the condition \rf{2.16} holds, and in this case it has two roots $h_1(a)\le h_*\le h_2(a)$. Moreover, whenever the inequality in \rf{2.16} is strict, i.e. $ a < a_{\max}(U_0)$, these  two roots are distinct.
\item[(iii)] According to the part (ii), whenever \rf{2.16} is fulfilled, there are two states  $U_i(a)=(h_i(a),u_i(a),a)$, where $u_i(a)=h_0u_0/h_i(a), i=1,2$ to which a stationary contact from  $U_0$ is possible. Moreover, the locations of these states can be determined as follows
$$
\begin{array}{ll}&U_1(a)\in G_1\quad\textrm{if}\quad u_0>0,\\
&U_1(a)\in G_3\quad\textrm{if}\quad u_0<0,\\
&U_2(a)\in G_2.\\
\end{array}$$
\end{itemize}
\end{lemma}

\noindent{\it Proof.} The parts (i) and (ii) can be easily deduced from the above argument. To prove (iii), it is sufficient to show that along the projection of $\WW_3(U_0)$ on the $(h,u)$-plane, the point $U_{\min}(U_0)=(h_{\min}(U_0), u_{\min}(U_0):=h_0u_0/h_{\min}(U_0))$, where $h_{\min}(U_0)$ is defined by \rf{2.15}, belongs to $\CC_+$ if $u_0>0$ and belongs to $\CC_-$ if $u_0<0$, and that $U_i(a)\in \WW_3(U_0), i=1,2,$ such that $U_2(a)\in G_2$ and $U_1(a)$ is located on the other side of $U_2(a)$ with respect to $\CC$. Indeed, let us define a function taking values along the stationary curve $\WW_3(U_0)$:
$$
\sigma(h):= u(h)^2 - gh ={h_0^2u_0^2\over h^2}-gh.
$$
Clearly, a point $U=(h,u,a)$ belongs to $G_1\cup G_3$ if and only if $\sigma(h)>0$ and $U$
 belongs to $G_2$ if and only if $\sigma(h)<0$. Since $\sigma(h_{\min}(U_0))=0$, the point
$U_{\min}(U_0)$ belongs to $\CC$. Obviously, $U_{\min}(U_0)\in \CC_+$ if $u_0>0$,
and $U_{\min}(U_0)\in \CC_-$ if $u_0<0$.  Thus, it remains to check  that
\begin{equation}
\sigma(h_1(a))>0,\quad \sigma(h_2(a))<0.
 \label{2.17}
\end{equation}
Since
$$
\sigma(h_{\min}(U-0))=0,\quad \sigma'(h)={-2h_0^2u_0^2\over h^3}-g<0,
$$
we can see that \rf{2.17} holds if
\begin{equation}
h_1(a)<h_{\min}(U_0)<h_2(a).
 \label{2.18}
\end{equation}
On the  other hand, we have
\begin{equation}
\begin{array}{ll}&\varphi(h)>0,\quad \textrm{if } h< h_1(a) \quad\textrm{or}\quad h>h_2(a),\\
&\varphi(h)<0,\quad  \textrm{if } h_1(a)<h<h_2(a).\\
\end{array} \label{2.19}
\end{equation}
And we have
$$
\varphi(h_{\min}(U_0))= 3(h_0u_0)^2 + (2g(a-a_0-h_0)-u_0^2){(h_0u_0)^{4/3}\over g^{2/3}}.
$$
It is a straightforward calculation to show that the condition
$$
a<a_{\max}(U_0)
$$
is equivalent to $$\varphi(h_{\min}(U_0))<0.$$ This together with  \rf{2.19} establish \rf{2.18}. Lemma \ref{lem21} is completely proved.
\qed

From Lemma \ref{lem21}, we can construct two-parameter wave sets.
The Riemann problem may therefore admit up to a one-parameter family of solutions.
To select a unique solution, we impose an admissibility condition for stationary contacts,
referred to as the {\it Monotonicity Criterion} and defined as follows:
\begin{itemize}
 \item[(MC)] Along any
stationary curve $\WW_3(U_0)$, the bottom level $a$ is monotone as
a function of $h$.  The total variation of the bottom level
component of any Riemann solution must not exceed  $|a_L-a_R|$, where $a_L, a_R$ are left-hand and
right-hand bottom levels.
\end{itemize}

A similar criterion was used
\cite{IsaacsonTemple92, IsaacsonTemple95},
\cite{LeFlochThanh03.2}, and
\cite{GoatinLeFloch}.

\begin{lemma} \label{lem22} The Monotonicity Criterion implies that any stationary shock
does not cross the boundary of strict hyperbolicity, in other
words:
\begin{itemize}
 \item[(i)] If $U_0\in G_1\cup G_3$, then
only the stationary contact based on the value $h_1(a)$ is
allowed, and one sets $\bar h(a)=h_1(a)$.
 \item[(ii)] If $U_0\in G_2$, then only the stationary
contact using $h_2(a)$ is allowed, and one sets $\bar h(a)=h_2(a)$.
\end{itemize}
\end{lemma}

Thus, $\bar h(a)$ is the admissible $h$-value of a right-hand state $U=(h=\bar h(a),u,a)$ of the stationary wave from a given left-hand state $U_0=(h_0,u_0,a_0)$.

%-----------------------------------------------------------------------------

{\it Proof.}
Recall that the Rankine-Hugoniot relations associate the linearly
degenerate field \rf{2.11} implies that the component $a$ can be
expressed as a function of $h$:
$$
 a=a(h)=a_0+{-u^2+u_0^2\over 2g}-h+h_0,
$$
where
$$
u=u(h)={h_0u_0\over h}.
$$
Thus, taking the derivative of $a$ with respect to $h$, we have
$$
\begin{array}{ll}a'(h)&={-uu'(h)\over g}-1=u{h_0u_0\over gh^2}-1\\
&={u^2\over  gh}-1={(u^2-gh)\over gh}\\
\end{array}$$
which has the same sign as $u^2-gh$. Thus, $a=a(h)$ is increasing with respect to $h$ in the domains $G_1, G_3$ and is decreasing in the domain $G_2$.  Thus, in order that
$a=a(h)$ is monotone as a function of $h$, the point $(h,u,a)$ must stay in the closure of the domain containing  $(h_0,u_0,a_0)$. The
conclusions of (i) and (ii) then follow.
\qed

We now explain how to compute the roots of the equation \rf{2.13}.
The above argument shows that whenever \rf{2.16} is satisfied, the equation \rf{2.13} admits
two roots $h_1(a), h_2(a)$ satisfying
\begin{equation}
h_1(a)\le h_{\min} =\Big({h_0^2u_0^2\over g}\Big)^{1/3}\le h_*=\frac{u_0^2+2g(a_0+h_0-a)}{3g} \le h_2(a)
  \label{2.20}
\end{equation}
and the inequalities are all strict whenever the inequality in \rf{2.16} is strict. Since $0<h_1(a)\le h_{\min}\le h_*$ and
$$
\begin{array}{ll}
&\varphi(0)\ge 0,\\
&\varphi(h_*)\le 0, \, \varphi(h_{\min})\le 0,\\
\end{array}
$$
the root $h_1(a)$ of \rf{2.13} can be computed, for instance
using the regula falsi method with the starting interval   $[0,h_{\min}]$, or $[0,h_*]$.
And since $h_2(a)\ge h_*$ and $\varphi'(h)>0,\varphi"(h)>0, h>h_*$, the root $h_2(a)$ can be computed using Newton's method with any starting point larger than $h_*$. We summarize this in the following lemma.

\begin{proposition} [Water height of stationary contacts]
\label{lem23}
The root $h_1(a)$ of \rf{2.13} can be computed using the regula falsi method for the starting interval $[0,h_{\min}]$, where $h_{\min} =\Big({h_0^2u_0^2\over g}\Big)^{1/3}$, or $[0,h_*]$, where $h_*=\frac{u_0^2+2g(a_0+h_0-a)}{3g}$, while the root $h_2(a)$ can be computed using Newton's method with any starting point larger than $h_*$.
\end{proposition}

To conclude this section, we point out that certain physical applications may actually require 
a different jump relation for the nonconservative product ---especially allowing for energy dissipation. This issue will not be addressed further in the present paper, however.

%=======================

\section{The Riemann problem revisited}

From the general theory of nonconservative  systems of balance laws, it is
known that if Riemann data belongs to a sufficiently small ball in a strictly hyperbolic region,
 then the Riemann problem admits a unique solution. It is worth to note that this result no longer
holds if any of these assumptions fails, for instance due to resonance.

Our goal in this section to to provide all possible explicit constructions for Riemann solutions, investigating when data are around the strictly hyperbolic boundary $\CC_\pm$.
There are several improvements in the construction of Riemann solutions in this paper over the ones in our previous work \cite{LeFlochThanh07}. First, we can determine  larger domains of existence by combining constructions in \cite{LeFlochThanh07} together. 
Second, the domains  where there is a unique solution or there are several solutions are
precisely determined.
Under the transformation $x\mapsto -x, u\mapsto -u$, a left-hand state $U=(h,u,a)$ in
$G_2$ or $G_3$ will be transferred to the right-hand state $V=(h,-u,a)$ in
$G_2$ or $G_1$, respectively. Thus, the construction for Riemann data around $\CC_-$ can be obtained from the one for Riemann data around $\CC_+$. We thus construct only the case where Riemann data are in $G_1\cup \CC_+\cup G_2$ and we separate into two regimes on which a corresponding construction based on the left-hand state $U_L$ is given:

\begin{itemize}
\item[$\bullet$] {\it Regime (A):} $U_L\in G_1\cup \CC_+$;
\item[$\bullet$] {\it Regime (B):} $U_L\in G_2$;
\end{itemize}
For each construction, depending on the location of the right-hand states $U_R$ and the sign $a_R-a_L$ there will be different types of solutions or the results on the existence and uniqueness.

As in \cite{LeFlochThanh07}, to solve \rf{1.1}--\rf{1.3} we project all the wave curves on the $(h,u)$-plane.

\bigskip

\section*{Notations}

\begin{itemize}
\item[(i)]  $U^{0}$ denotes the state resulted from a stationary contact wave from $U$;
\item[(ii)] $U^{\#}$ is the state defined  in Lemma \ref{lem31} so that $\bar\ld_1(U,U^{\#})=0$;
\item[(iii)]
$W_k(U_i,U_j) $ ($S_k(U_i,U_j), R_k(U_i,U_j)$) denotes the $k$th-wave ($k$th-shock, $k$th-rarefaction wave, respectively) connecting the left-hand state $U_i$ to the right-hand state $U_j$, $k=1,2,3$;
\item[(iv)]
$W_m(U_i,U_j)\oplus W_n(U_j,U_k)$ indicates that there is an $m$th-wave from the left-hand state $U_i$ to the right-hand state $U_j$, followed by an $n$th-wave from the left-hand state $U_j$ to the right-hand state $U_k$, $m,n\in\{1,2,3\}$.
\end{itemize}

%---------------------------------------------------------------

\subsection{Regime (A). Eigenvalues at $U_L$ with coinciding signs}

Let $\bar G_1$ denote the closure of $G_1$. We assume that  $U_L\in \bar G_1$, or equivalently $\ld_i(U_L)\ge 0, i=1,2,3$.

\noindent\underline{\it Construction (A1).} In this case (the projection on $(h,u)$-plane of) $U_R$ is located in a "higher" region containing $U_L$ in the $(h,u)$-plane.

\begin{figure}
  \includegraphics[width=10truecm]{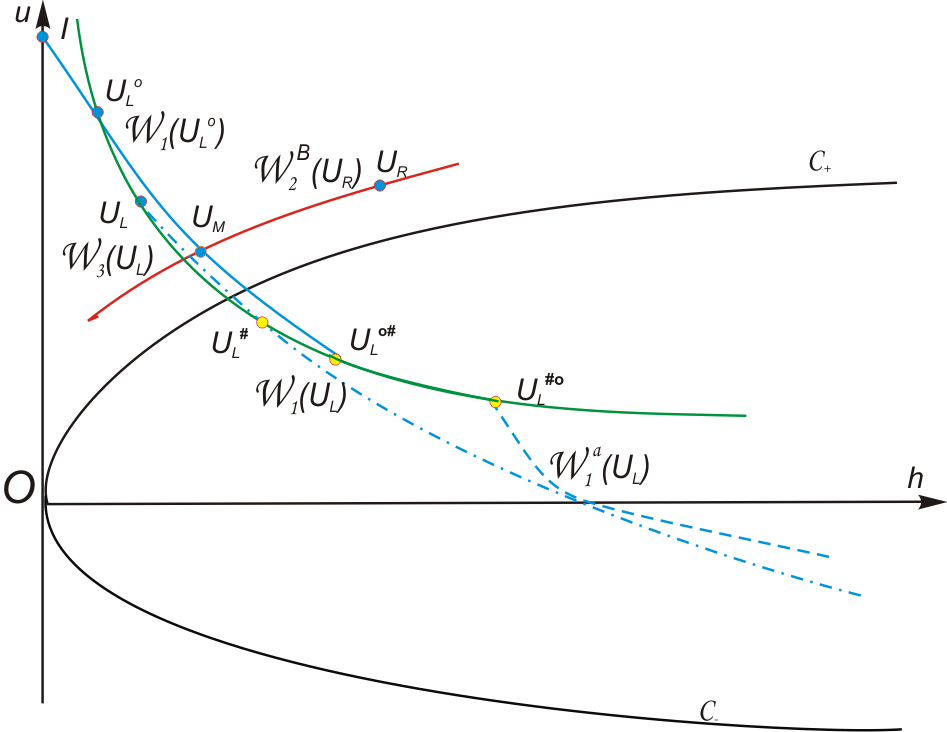}\\
  \caption{Construction (A1), $a_L>a_R$: a solution of the form \rf{3.1}.}\label{Riem1a}
\end{figure}

If $a_L\ge a_R$ (or $a_L< a_R\le a_{\max}(U_L)$), the solution begins with a stationary contact upward (downward, respectively) along $\WW_3(U_L)$  from $U_L$ to the state $U_L^o\in \WW_3(U_L)\cap G_1$, shifting the level $a_L$ directly to the level $a_R$.
Let
$$
\{U_M=(h_M,u_M,a_R)\}=\WW_1(U_L^o)\cap \WW_2^B(U_R).
$$
Providing that $\bar\ld_1(U_L^o,U_M)\ge 0$, or equivalently, as seen from Lemma \ref{lem31}, $h_M\le h_L^{o\#}$,
 the solution can continue by a $1$-wave from $U_L^o$ to $U_M$, followed by a $2$-wave from $U_M$ to $U_R$. Thus, the solution is
\begin{equation}
W_3(U_L,U_L^o)\oplus W_1(U_L^o,U_M)\oplus W_2(U_M,U_R).
\label{3.1}
\end{equation}
See Figure \ref{Riem1a}.
This construction can be extended if $\WW_2^B(U_R)$ lies entirely above $\WW_1(U_L^o)$. In this case let
$I$ and $J$ be the intersection points of   $\WW_1(U_L^o)$ and  $\WW_2^B(U_R)$ with the axis $\{h=0\}$, respectively:
\begin{equation}
\{I\}= \WW_1(U_L^o)\cap \{h=0\},\quad \{J\}= \WW_2^B(U_R)\cap\{h=0\},
\label{3.2}
\end{equation}
 then the solution can be seen as a dry part $W_o(I,J)$ between $I$ and $J$. Thus, the solution in this case is
$$
W_3(U_L,U_L^o)\oplus R_1(U_L^o,I)\oplus W_o(I,J)\oplus R_2(J,U_R).
$$

\begin{rmk} As seen by Lemma \ref{lem21}, if $a_L<a_R$, the condition
$$
a_R\le a_{\max}(U_L)
$$
is necessary for the stationary contact $W_3(U_L,U_L^o)$. Therefore, if this condition fails, there is no solution even if $U_L=U_R$. The necessary and sufficient conditions for the existence of the solution \rf{3.1} is that $U_L^{o\#}$ is located below or on the curve $\WW_2^B(U_R)$. This domain clearly covers a large crossing-strictly-hyperbolic-boundary neighborhood of $U_L$.
\end{rmk}

\begin{figure}
  \includegraphics[width=10truecm]{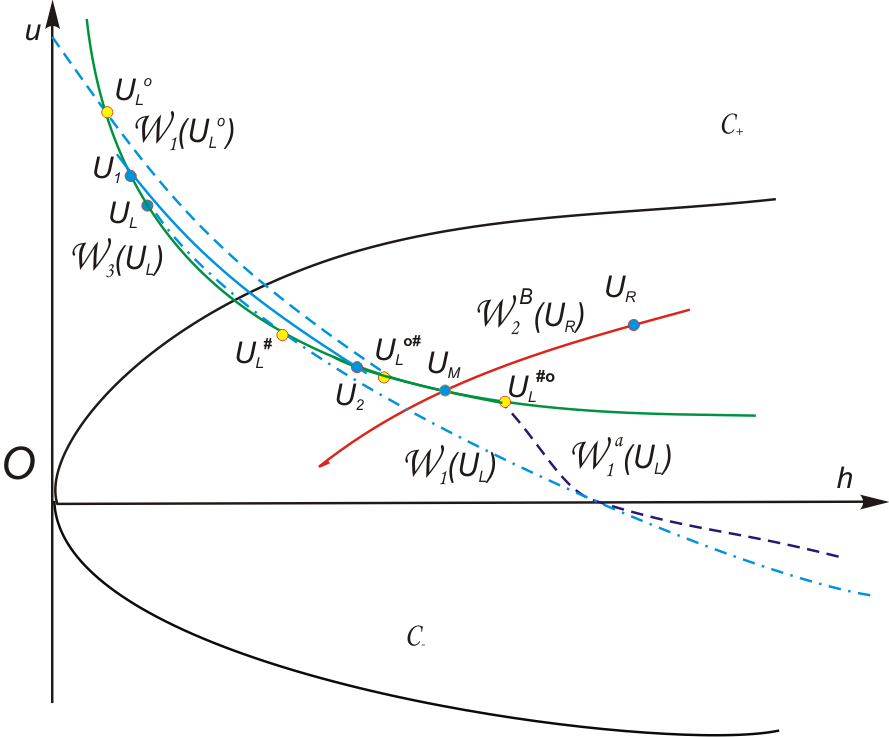}\\
  \caption{Construction (A2), $a_L>a_R$: a solution of the form \rf{3.3}.}\label{Riem1b}
\end{figure}

 \noindent\underline{\it Construction (A2).}
In this construction, we will see an interesting phenomenon when  wave speeds associate with different characteristic fields coincide.
 Roughly speaking, this case concerns with the fact that $U_R$ moves limitedly downward from the case $G_1$. Instead of using
``complete" stationary contact from $U_L$ to $U_L^o$ as in the first possibility, the solution now begins with a
``half-way" stationary contact $W_3(U_L,U_1)$ from $U_L=(h,u,a_L)$ to some state $U_1=U_L^o(a)=(h,u,a)\in \WW_3(U_L)$, where $a$ between $a_L$ and $a_R$. The solution then continues by a $1$-shock wave with zero speed from $U_1$ to $U_2=U_1^{\#}\in \WW_1(U_1)\cap G_2$. Observe that $U_2$ still belongs to $\WW_3(U_L)$, since $h_2u_2=h_1u_1=h_Lu_L$, as indicated by Lemma \ref{lem31}. The solution continues by a stationary contact from $U_2$ to a state $U_M(a)\in \WW_3(U_L)$.
The set of these points $U_M(a), a\in [a_L,a_R]$ forms a curve pattern denoted by $\LL$.
Whenever
$$
\WW_2^B(U_R)\cap \LL\ne\emptyset
$$
there is a solution containing three discontinuities having the same zero speed of the form
\begin{equation}
W_3(U_L,U_1)\oplus S_1(U_1,U_2)\oplus W_3(U_2 ,U_M)\oplus W_2(U_M,U_R).
\label{3.3}
\end{equation}
See Figure \ref{Riem1b}.

 \begin{rmk}
 The necessary and sufficient conditions for the existence of the solution \rf{3.3} is that $U_L^{o\#}$ is located above or on the curve $\WW_2^B(U_R)$, and $U_L^{\# o}$ is located below or on the curve $\WW_2^B(U_R)$. This domain  covers a region in $G_2^+$ and $G_1$ which is
``far away" from $U_L$.
\end{rmk}

 It is interesting that at the limit $a=a_R$ at the first jump, we get the first possibility. If $a=a_L$, then the solution simply begins with a $1$-shock wave with zero speed followed by a stationary contact  shifting $a$ from $a_L$ to $a_R$. This limit case can be connected to the following possibility.

\begin{figure}
  \includegraphics[width=10truecm]{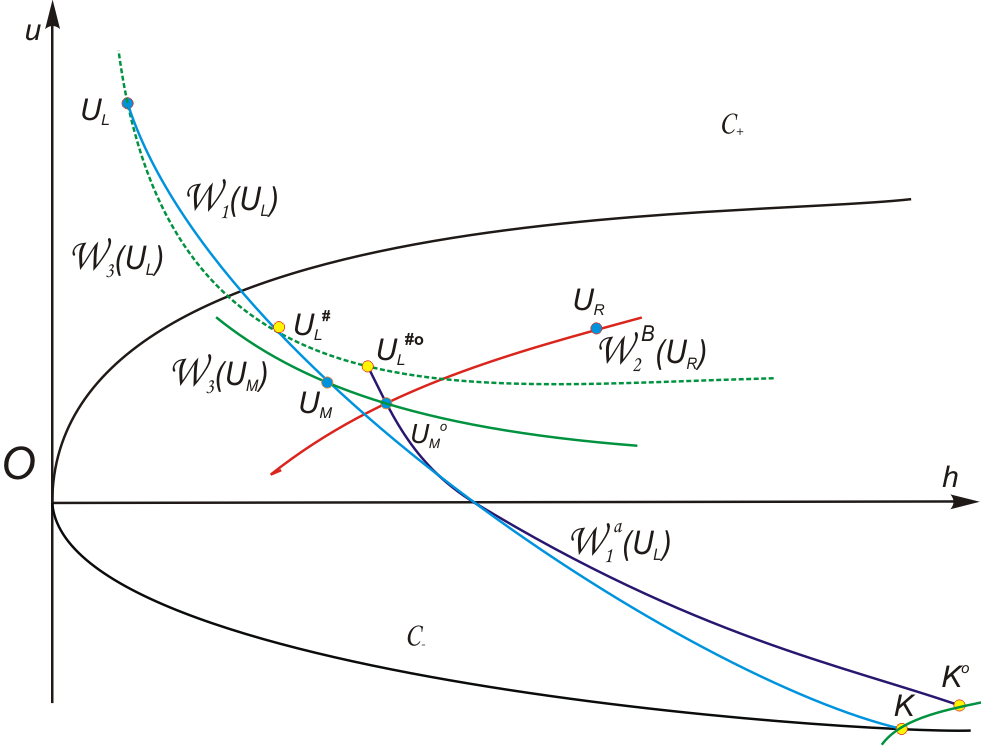}\\
  \caption{Construction (A3), $a_L>a_R$: a solution of the form \rf{3.5}.}\label{Riem1c}
\end{figure}

\noindent\underline{\it Construction (A3).}
The solution begins with  a  strong $1$-shock wave from $U_L$ to any state $U\in \WW_1(U_L)\cap G_2$ such that $\bar\ld_1(U_L,U)\le 0$. This shock wave is followed by a stationary contact to a state $U^o$ shifting $a$ from $a_L$ to $a_R$. The set of these states $U^o$ form a curve denoted by $\WW_1^a(U_L)$. That is
\begin{equation}
\begin{array}{ll}\WW_1^a(U_L):=\{U^o: &\exists W_3(U, U^o)\quad\textrm{shifting } a_L \quad \textrm{to } a_R,\\
 &U=(h,u,a_L)\in \WW_1(U_L), \ld_1(U)\le 0\}.\\
 \end{array}
\label{3.14}
\end{equation}
Whenever
\begin{equation}
\emptyset\ne \WW_2^B(U_R)\cap \WW_1^a(U_L)=\{U_M^o\}\subset G_2,\quad\textrm{and}\quad \bar\ld_2(U_M^o,U_R)\ge 0,
\label{3.4}
\end{equation}
there will be a Riemann solution defined by
\begin{equation}
S_1(U_L,U_M)\oplus W_3(U_M,U_M^o)\oplus W_2(U_M^o,U_R).
\label{3.5}
\end{equation}
See Figure \ref{Riem1c}.
In the limit case of \rf{3.3} where $U_1\equiv U_L$,  the solution \rf{3.3} coincides with the solution \rf{3.5}.

Let $K$ denote the lower limit state on $\WW_1(U_L)$ that the solution \rf{3.5} makes sense, and let $K^o\in G_2$ denote   the right-hand state resulted from  a stationary contact from $K$ shifting $a_L$ to $a_R$. Thus, we have
\begin{equation}
\begin{array}{ll}&\{K\}=\WW_1(U_L)\cap \CC_-,\quad \textrm{if } a_L\ge a_R,\\
&K\in \WW_1(U_L)\quad\textrm{such that } a_{\max}(K)=a_R,\quad \textrm{if } a_L< a_R.\\
\end{array}\label{3.6}
\end{equation}

 \begin{rmk}
 The solution \rf{3.5} makes sense if $U_L^{\# o}$ is above or on the curve $\WW_2^B(U_R)$, and  $K^o$ lies below or on the curve $\WW_2^B(U_R)$ and $\bar\ld_3(K^o,U_R)\ge 0$.
\end{rmk}

The union of the wave patterns $\WW_1(U_1)\cup \LL\cup \WW_1^a(U_L)$ form a continuous curve.   The Riemann problem thus admits a solution whenever $\WW_2^B(U_R)$ intersects  $\WW_1(U_1)\cup \LL\cup \WW_1^a(U_L)$ or  $\WW_2^B(U_R)$ intersects  with $\{h=0\}$ at a point above the point $I$.
We can see that this happens for a large domain of $U_R$ containing $U_L$.
See Figures \ref{Riem1a}, \ref{Riem1b}, and \ref{Riem1c}.

\begin{figure}
  \includegraphics[width=10truecm]{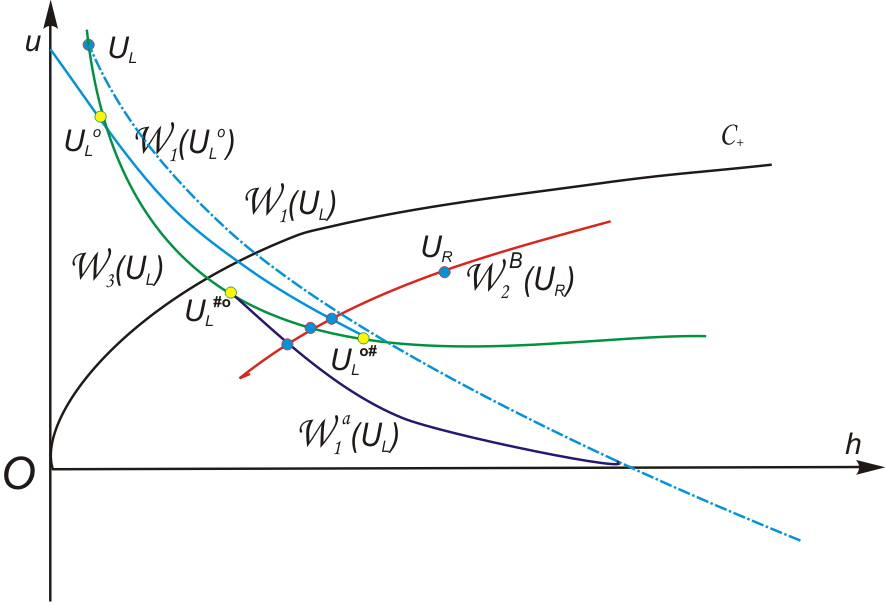}\\
  \caption{Non-uniqueness:
three admissible Riemann solutions of the form \rf{3.1}, \rf{3.3}, and \rf{3.5}.}\label{Riem1-multiple}
\end{figure}

If the wave pattern $\LL$ lies entirely on one side with respect to the curve $\WW_2^B(U_R)$, then $\WW_2^B(U_R)$ intersects either $\WW_1(U_1)$ or $\WW_1^a(U_L)$ at most one point. Therefore, then \rf{3.1} or \rf{3.5} is the {\it unique solution}. Besides, if $\WW_2^B(U_R)$ intersects the wave pattern $\LL$, and  if $h_L^{\#o}\ge  h_L^{o\#}$,
then the point  $U_L^{\#o}$ is located below the point $U_L^{o\#}$ on the curve $\WW_3(U_L)$. Thus, the curve $\WW_2^B(U_R)$ does not meet  $\WW_1(U_1)$ nor $\WW_1^a(U_L)$, except possibly at the endpoints
  $U_L^{\#o}\in \LL$ and $U_L^{o\#}\in\LL$. In this case, \rf{3.3} is the sole solution.
  In summary, the Riemann problem for \rf{1.1}--\rf{1.2} always {\it has at most one solution} whenever  $h_L^{\#o}\ge  h_L^{o\#}$.

In the case where $h_L^{\#o}< h_L^{o\#}$, {\it there can be three solutions}, as $\WW_2^B(U_R)$ can meet all the three curve patterns $\WW_1(U_1),  \LL$ and $\WW_1^a(U_L)$, or
\begin{equation}
h_L^{\#o}< h_L^{o\#},\quad \Phi_2(U_L^{\#o};U_R)>0>\Phi_2(U_L^{o\#};U_R),
\label{3.7}
\end{equation}
where the function $\Phi_2(U;U_R)$ is defined by \rf{2.10}.
See Figure \ref{Riem1-multiple}.

The above argument leads us to the following theorem.

\begin{theorem}[Riemann problem for the shallow water equations]
 \label{theo32} Given a left-hand state $U_L\in G_1$. Depending on the location of the right-hand state $U_R$ we have the following conclusions.
\begin{itemize}
\item[(a)] {\bf Existence.} The Riemann problem \rf{1.1}--\rf{1.3} admits a solution if $Q^o$ defined in Construction (A3) lies below or on the curve $\WW_2^B(U_R)$, and that if $\WW_2^B(U_R)$ intersects with $\WW_1^a(U_L)$ at some point $U_M^o\in G_2^-$, then   $\bar\ld_2(U_M^o,U_R)\ge 0$.

\item[(b)] {\bf Regime of uniqueness.} The Riemann problem \rf{1.1}--\rf{1.3} has at most one solution if

--  either $h_L^{\#o}\ge  h_L^{o\#}$;

-- or  $h_L^{\#o}<  h_L^{o\#}$, and the states  $U_L^{\#o}$ and $U_L^{o\#}$ are located on the same side with respect to the curve $\WW_2^B(U_R)$.

\item[(c)] {\bf Multiple  solutions.} If $h_L^{\#o}<  h_L^{o\#}$, and if the state  $U_L^{\#o}$ lies above the curve $\WW_2^B(U_R)$  while the state  $U_L^{o\#}$ lies below the curve $\WW_2^B(U_R)$,  then the Riemann problem \rf{1.1}--\rf{1.3} has three solutions.
\end{itemize}
\end{theorem}

\noindent{\it Example.} We provide some numerical experiments to illustrate two situations: $h_L^{\#o} >  h_L^{o\#}$, and $h_L^{\#o} <  h_L^{o\#}$ corresponding to the two cases $a_L> a_R$ (see Tables A1-A3) and $a_L< a_R$ (see Tables A4, A5). We take at random the state $U_L$ and $a_R$.

(a) $a_L> a_R$: all experiments show that  $h_L^{\#o}> h_L^{o\#}$.
 In Table A1, $U_L\in G_1$.
$$
\framebox{
\begin{tabular}{l l c l c l c}
{\bf Table A1} \\
\hline
   States          &\vline  &  $U_L$   &\vline &   $U_L^{\#o}$ & \vline &   $U_L^{o\#}$ \\
\hline
Water Height $h$   & \vline &     0.5  & \vline &   1.1930011   &\vline &   1.1171275   \\
Velocity  $u$      & \vline &     4    & \vline &   1.6764444   &\vline &   1.790306    \\
Bottom Level  $a$  & \vline &     1    & \vline &   0.9         &\vline &   0.9         \\
\end{tabular}
}
\label{TableA1}
$$
In Table A2, $U_L\in \CC_+$.
$$
\framebox{
\begin{tabular}{l l c l c l c}
{\bf Table A2} \\
\hline
   States          &\vline  &  $U_L$         &\vline &   $U_L^{\#o}$ & \vline &   $U_L^{o\#}$ \\
\hline
Water Height $h$   & \vline &     1         & \vline &   1.3075478  &\vline &   1.2558035  \\
Velocity  $u$      & \vline &     3.1304952  & \vline &  2.3941726   &\vline &  2.4928225   \\
Bottom Level  $a$  & \vline &     1         & \vline &   0.9  &\vline &   0.9  \\
\end{tabular}
}
\label{TableA2}
$$
In Table A3, $U_L\in G_1$ is far away from $\CC_+$.
$$
\framebox{
\begin{tabular}{l l c l c l c}
{\bf Table A3} \\
\hline
   States          &\vline  &  $U_L$         &\vline &   $U_L^{\#o}$ & \vline &   $U_L^{o\#}$ \\
\hline
Water Height $h$   & \vline &     0.01  & \vline & 0.54763636  &\vline & 0.44902891   \\
Velocity  $u$      & \vline &     10    & \vline & 0.18260292  &\vline & 0.22270281   \\
Bottom Level  $a$  & \vline &     1     & \vline & 0.9  &\vline & 0.9 \\
\end{tabular}
}
\label{TableA3}
$$

(b) $a_L< a_R$: all experiments show that  $h_L^{\#o}< h_L^{o\#}$.
In Table A4, $U_L\in G_1$.
$$
\framebox{
\begin{tabular}{l l c l c l c}
{\bf Table A4} \\
\hline
   States          &\vline  &  $U_L$   &\vline &   $U_L^{\#o}$ & \vline &   $U_L^{o\#}$ \\
\hline
Water Height $h$   & \vline &     0.5  & \vline &   0.86127059  &\vline &   0.96534766   \\
Velocity  $u$      & \vline &     4    & \vline &   2.3221506  &\vline &   2.0717925   \\
Bottom Level  $a$  & \vline &     0.9    & \vline & 1    &\vline & 1     \\
\end{tabular}
}
\label{TableA4}
$$
In Table A5, $U_L\in G_1$ is far away from $\CC_+$.
$$
\framebox{
\begin{tabular}{l l c l c l c}
{\bf Table A5} \\
\hline
    States          &\vline  &  $U_L$         &\vline &   $U_L^{\#o}$ & \vline &   $U_L^{o\#}$ \\
\hline
Water Height $h$   & \vline &     0.01  & \vline &  1.2748668 &\vline  & 1.3718425   \\
Velocity  $u$      & \vline &     10    & \vline &  0.78439566 &\vline  & 0.72894668   \\
Bottom Level  $a$  & \vline &     0.9    & \vline & 1 &\vline &1\\
\end{tabular}
}
\label{TableA5}
$$

\begin{rmk} We conjecture that if $a_L> a_R$, then $h_L^{\#o}> h_L^{o\#}$, and if  $a_L< a_R$, then $h_L^{\#o}< h_L^{o\#}$.
If this conjecture holds, then Theorem \ref{theo32} implies that when $a_L\ge a_R$, the Riemann problem 
has at most one solution for $U_L\in G_1$.
\end{rmk}

\subsection{Regime (B). Eigenvalues at $U_L$ with opposite signs}

In this subsection we consider the case where the left-hand state $U_L$ moves downward from the Regime (A): $U_L\in \bar G_2$, or $\ld_1(U_L)<0=\ld_3(U_L)<\ld_2(U_L)$.

\noindent\underline{\it Construction (B1).}   For  $U_R$ in a ``higher'' position, there can be two types of solutions depending on whether $a_L\ge a_R$.

\begin{figure}
  \includegraphics[width=10truecm]{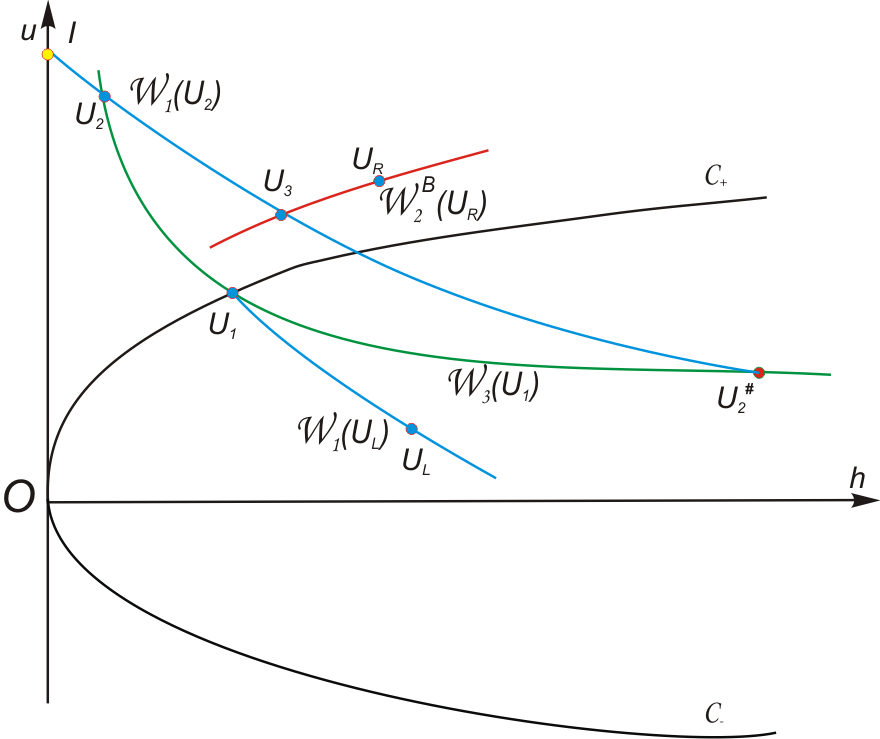}\\
  \caption{Construction (B1), $a_L>a_R$: a solution of the form \rf{3.8}.}\label{Riem2a}
\end{figure}

If $a_L> a_R$ a solution can be constructed as follows.
The  solution begins from $U_L$ with a 1-rarefaction wave until it reaches $\CC_+$ at a state $U_1\in \WW_1(U_L)\cap \CC_+$. A straightforward calculation gives
 $$
 U_1=\Big(\big({u_L\over 3\sqrt{g}}+{2\over 3}\sqrt{h_L}\big)^2, {1\over 3}u_L+{2\over 3}\sqrt{gh_L},a_L\Big).
 $$
  This rarefaction wave can be followed by a stationary jump $W_3(U_1,U_2)$ into $G_1$. This stationary wave is possible since $a_L\ge a_R$.  Let $\{U_3\}=\WW_1(U_2)\cap \WW_2^B(U_R)$. The solution is then continued by a 1-wave from $U_2$ to $U_3$, followed by a 2-wave from $U_3$ to $U_R$.
Thus, the solution is given by the formula
\begin{equation}
R_1(U_L,U_1)\oplus W_3(U_1,U_2) \oplus W_1(U_2, U_3)\oplus W_2(U_3,U_R).
\label{3.8}
\end{equation}
See Figure \ref{Riem2a}.
The construction makes sense if $\bar\ld_1(U_2,U_3)\ge 0$, which means $U_3$ has to be above $U_2^{\#}$ on $\WW_1(U_2)$.  This construction can also be extended if $\WW_2^B(U_R)$ lies entirely above $\WW_1(U_2)$. In this case let
$I$ and $J$ be the intersection points of   $\WW_1(U_2)$ and  $\WW_2^B(U_R)$ with the axis $\{h=0\}$, respectively:
$$%
\{I\}= \WW_1(U_2)\cap \{h=0\},\quad \{J\}= \WW_2^B(U_R)\cap\{h=0\}.
$$%
 Then, the solution can be seen as containing a dry part $W_o(I,J)$ between $I$ and $J$. Thus, the solution in this case is
\begin{equation}
R_1(U_L,U_1)\oplus W_3(U_1,U_2) \oplus W_1(U_2,I)\oplus W_o(I,J)\oplus R_2(J,U_R).
\label{3.9}
\end{equation}

\begin{figure}
  \includegraphics[width=10truecm]{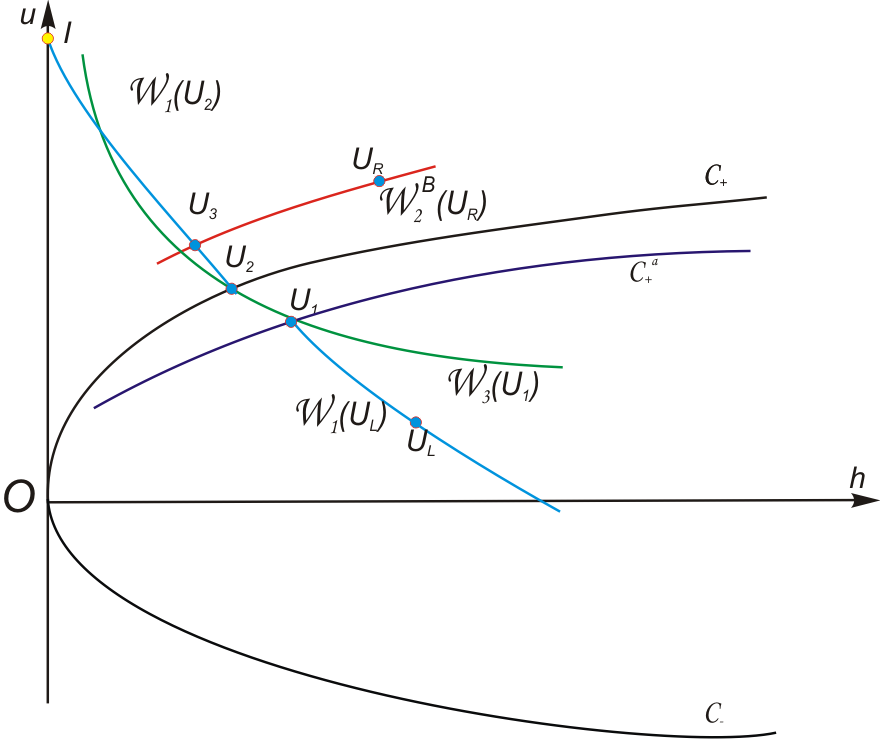}\\
  \caption{Construction (B1), $a_L\le a_R$: a solution of the form \rf{3.10}.}\label{Riem2ab}
\end{figure}

If $a_L\le a_R$ a solution of another type can be constructed as follows.
To each $U\in \CC_+$, a stationary contact to $U^o\in G_2$ downing $a=a_R$  to $a=a_L$ is possible, since $a_R>a_L$. The set of all these states $U^o$ form a curve denoted by $\CC_+^a$.
Let
$$
\{U_1\}=\WW_1(U_L)\cap \CC_+^a.
$$
Then, the  solution begins by a $1$-wave $W_1(U_L,U_1)$,  followed by a stationary jump $W_3(U_1,U_2=U_1^o)$ to $U_2\in\CC_+$. Let $\{U_3\}=\WW_1(U_2)\cap \WW_2^B(U_R)$. The solution is then continued by a 1-rarefaction wave from $U_2$ to $U_3$, followed by a 2-wave from $U_3$ to $U_R$.
Thus, the solution is given by the formula
\begin{equation}
W_1(U_L,U_1)\oplus W_3(U_1,U_2) \oplus R_1(U_2, U_3)\oplus W_2(U_3,U_R).
\label{3.10}
\end{equation}
See Figure \ref{Riem2ab}.
The construction makes sense if $\ld_1(U_3)\ge 0$, or $U_3\in\bar G_1$.
This construction can also be extended if $\WW_2^B(U_R)$ lies entirely above $\WW_1(U_2)$. In this case let
$I$ and $J$ be the intersection points of   $\WW_1(U_2)$ and  $\WW_2^B(U_R)$ with the axis $\{h=0\}$, respectively:
$$
\{I\}= \WW_1(U_2)\cap \{h=0\},\quad \{J\}= \WW_2^B(U_R)\cap\{h=0\}.
$$
 Then, the solution can be seen as containing a dry part $W_o(I,J)$ between $I$ and $J$. Thus, the solution in this case is
\begin{equation}
W_1(U_L,U_1)\oplus W_3(U_1,U_2)\oplus R_1(U_2,I) \oplus W_o(I,J)\oplus R_2(J,U_R).
\label{3.11}
\end{equation}

The wave structure of the solutions \rf{3.8} and \rf{3.10} are the same, but the state at which the solution reaches the strictly hyperbolic boundary $\CC$ using a different wave. However, one may argue that in both cases the solution uses a stationary contact to reach $\CC_+$ from either side of $\CC_+$. Moreover, all the states in the solution $(U_L,U_1,U_2,U_3,U_R)$ can be in an arbitrarily small ball center on $\CC_+$.

\noindent\underline{\it Construction (B2).} This case holds only when $a_L>a_R$. Again, there is an interesting phenomenon as wave speeds associate with different characteristic fields coincide  and all equal zero. The solution therefore contain three waves with the same zero speed.

\begin{figure}
  \includegraphics[width=10truecm]{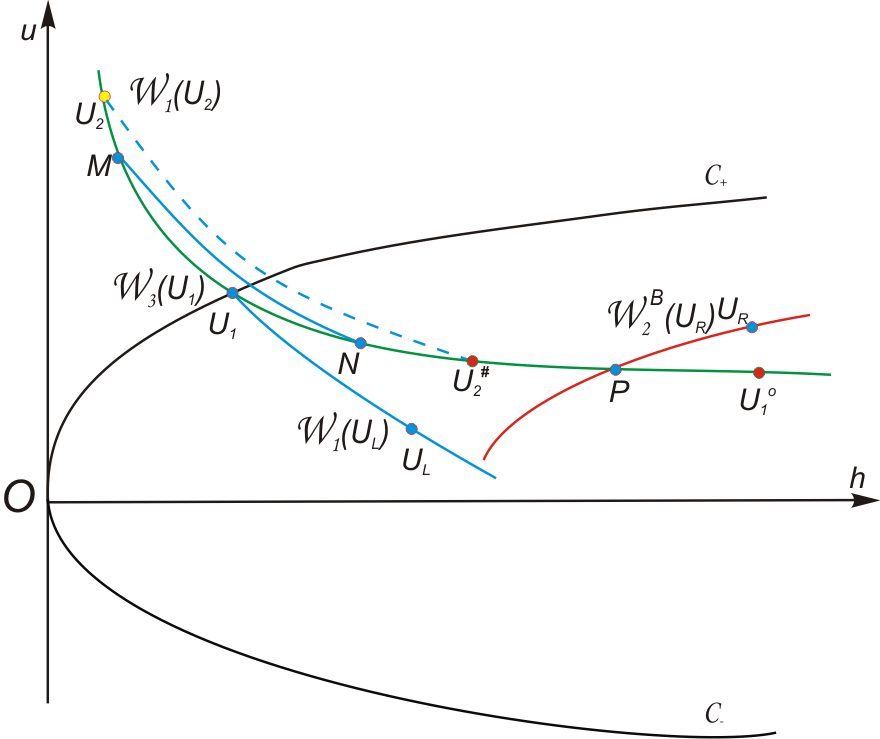}\\
  \caption{Construction (B2), $a_L>a_R$: a solution of the form \rf{3.12}.}\label{Riem2b}
\end{figure}

The solution begins with a $1$-rarefaction wave until it reached $\CC_+$ at $U_1$. At $U_1$, the solution may jump to $G_1$ using a
``half-way" stationary wave to a state $M=M(a)=U_1^{o}(a)$ from the bottom height $a_L$ to any $a\in [a_R,a_L]$. Then, the solution can continue by a $1$-shock with zero speed from $M$ to $N=N(a)=N^{\#}(a)\in G_2^+$, followed by a stationary wave from $N$ to $P=P(a)=N^o(a)$ with a shift in $a$-component from $a$ to $a_R$. The set of these states $P(a)$ form a curve pattern $\LL$. So, whenever $\emptyset\ne\WW_2^B(U_R)\cap \LL =\{P=P(a)\}$, there is a Riemann solution containing three zero-speed waves of the form
\begin{equation}
R_1(U_L,U_1)\oplus W_3(U_1,M(a))\oplus S_1(M(a),N(a))\oplus W_3(N(a), P(a))\oplus W_2(P(a),U_R).
\label{3.12}
\end{equation}
See Figure \ref{Riem2b}.

Observe that this solution coincides with the one in 
Construction (B1) if the first stationary wave from $U_1$ to $M=U_2$ shifts $a$-component from $a_L$ directly to $a_R$. The other limit case where the first stationary wave to $G_1$ is not used gives a connection to the following possibility.

\begin{figure}
  \includegraphics[width=10truecm]{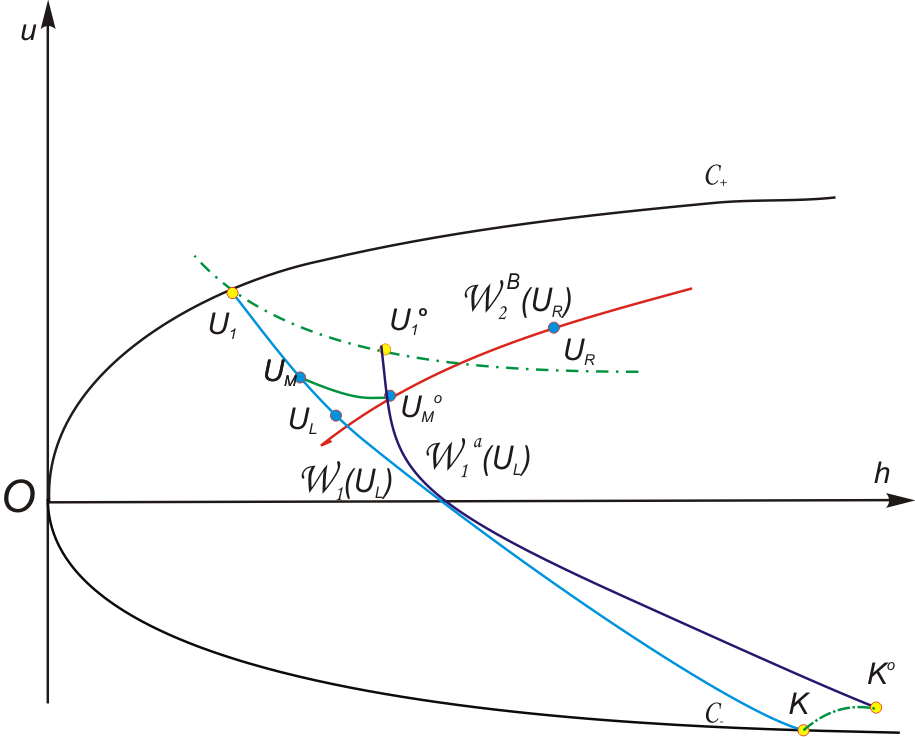}\\
  \caption{Construction (B3), $a_L>a_R$: a solution of the form \rf{3.15}.}\label{Riem2c}
\end{figure}

\noindent\underline{\it Construction (B3).} In this case the Riemann data can be altogether in a arbitrarily small ball in $G_2$. Assume first that $a_L\ge a_R$.
Let
\begin{equation}
\begin{array}{ll}&U_1=\WW_1(U_L)\cap \CC_+,\quad \textrm{and}\quad U_1^0\in G_2+ \quad \textrm{ resulted by } \quad W_3(U_1,U_1^o),\\
&K=\WW_1(U_L)\cap \CC_-,\quad \textrm{and}\quad K^0\in G_2- \quad \textrm{ resulted by} \quad W_3(K,K^o).\\
\label{3.13}
\end{array}\end{equation}
From any state $U\in \WW_1(U_L)$, where $\ld_1(U_L)\le 0$ ($U_L$ is below $U_1$ or coincides with $U_1$), we use a stationary jump to  a state $U^o$, shifting the bottom height from $a_L$ down to $a_R$. The set  of these states $U^o$ form a ``composite" curve $ \WW_1^a(U_L)$ as defined by \rf{3.14}.
The curve $\WW_1^a(U_L)$ is thus a path between $U_1^o$ and $K^o$.
Whenever $\emptyset\ne \WW_2^B(U_R)\cap \WW_1^a(U_L)=\{U_M^o\}$,  a Riemann solution can be determined by
\begin{equation}
W_1(U_L,U_M)\oplus W_3(U_M,U_M^o)\oplus W_2(U_M^o,U_R),
\label{3.15}
\end{equation}
where $U_M\in \WW_1(U_L)$, provided $U_R\in G_2$ or $\bar\ld_2(U_M^o,U_R)\ge 0$.
See Figure \ref{Riem2c}.

Second, consider the case $a_L< a_R$. Let $\CC_+^a$ as in the case for the solution of the type \rf{3.15}.
To each $U\in \CC_\pm$, a stationary contact to $U^o\in G_2$ downing back $a=a_R$  to $a=a_L$ is possible, since $a_R>a_L$. The set of all these states $U^o$ form two curves denoted by $\CC_\pm^a$.
Let
\begin{equation}
\begin{array}{ll}&\{U_1\}=\WW_1(U_L)\cap \CC_+^a, \quad \textrm{and}\quad U_1^0\in \CC_+ \quad\textrm{resulted by} \quad W_3(U_1^o,U_1),\\
&K\in \WW_1(U_L), \quad \textrm{and}\quad K^0\in \CC_-\quad \textrm{resulted by}\quad  W_3(K^o,K)
\end{array}\label{3.16}
\end{equation}
decreasing $a_R$ to $a_L$.
 From any state $U\in \WW_1(U_L), \ld_1(U)\le \ld_1(U_1)$, there is a stationary jump to  a state $U^o$, shifting the bottom height from $a_L$  to $a_R$. The set  of these states $U^o$ form a  composite curve also denoted by $\WW_1^a(U_L)$.
Whenever $\emptyset\ne \WW_2^B(U_R)\cap \WW_1^a(U_L)=\{U_M^o\}$,  a Riemann solution can be determined by
\begin{equation}
W_1(U_L,U_M)\oplus W_3(U_M,U_M^o)\oplus W_2(U_M^o,U_R),
\label{3.17}
\end{equation}
where $U_M\in \WW_1(U_L)$, provided $U_R\in G_2$ or $\bar\ld_2(U_M^o,U_R)\ge 0$.

\begin{figure}
  \includegraphics[width=10truecm]{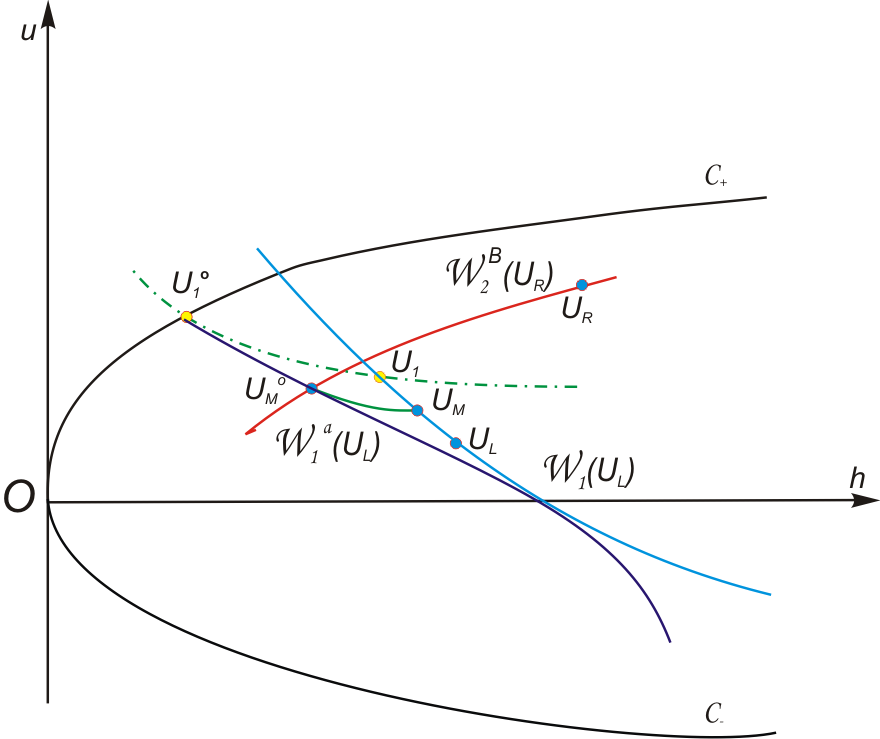}\\
  \caption{Construction (B3), $a_L<a_R$: a solution of the form \rf{3.17}.}\label{Riem2c2}
\end{figure}

\begin{rmk}
In both cases $a_L> a_R$ and $a_L\le a_R$, the condition for $\WW_2^B(U_R)\cap \WW_1^a(U_L)\ne \emptyset$ is that $U_1^0$ lies above $\WW_2^B(U_R)$ and $K^o$ lies below $\WW_2^B(U_R)$.
See Figure \ref{Riem2c2}.
\end{rmk}

Let us now  discuss the existence and uniqueness.
Assume first that  $a_L\le a_R$. In this case, only Constructions (B1) and (B3) are available. The limit case of \rf{3.8} of (B1) when $U_3\equiv U_2$ coincides with the limit case of \rf{3.15} of (B3). Thus, the union $\WW_1(U_2)\cup \WW_1^a(U_L)$ form a continuous decreasing curve (the curve can be considered as the graph of $u$ being a decreasing function of $h$) and that $\WW_1(U_2)$ and $\WW_1^a(U_L)$ meets only at one point $U_2$. Since $\WW_2^B(U_R)$ is an increasing curve, there always a unique intersection point of $\WW_2^B(U_R)$ and the union $\WW_1(U_2)\cup \WW_1^a(U_L)$ if $K^o$ lies below or on the curve $\WW_2^B(U_R)$. This implies that the Riemann problem for \rf{1.1}--\rf{1.2} always admits a unique solution if $K^o$ lies below or on the curve $\WW_2^B(U_R)$.

\begin{figure}
  \includegraphics[width=10truecm]{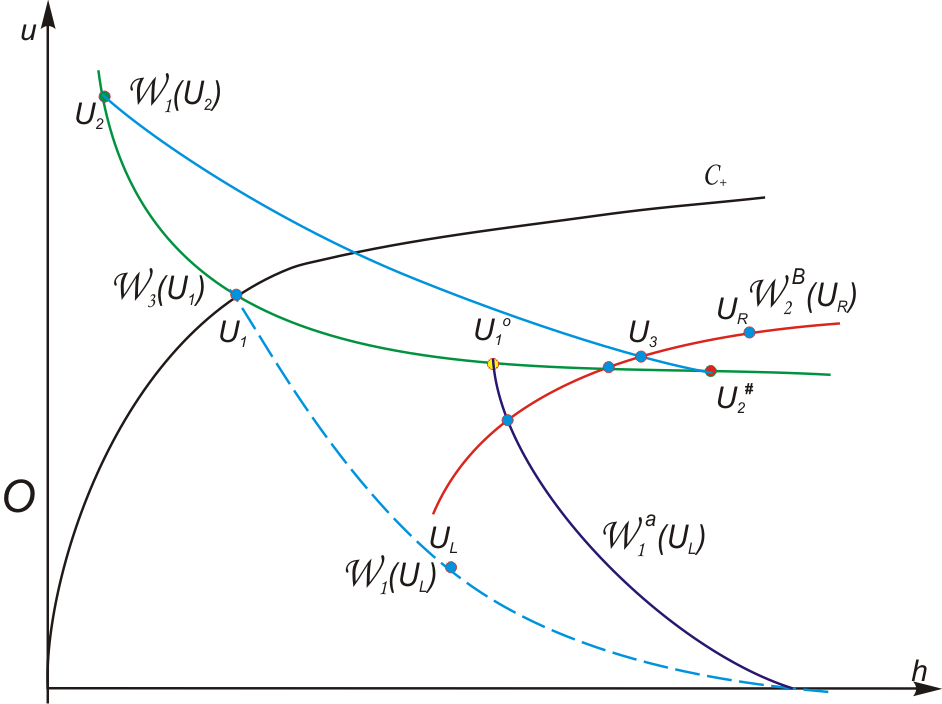}\\
  \caption{Regime (B):  whether $h_1^o< h_2^{\#}$
determines that the Riemann solution is unique 
or else there exist multiple solutions
$U_L\in G_2$ and $a_L>a_R$ (cf.~Theorem \ref{theo33}).}\label{Riem2-multiple}
\end{figure}

 Next, assume that $a_L>a_R$. Let $U_1^o\in G_2$ denote the state resulted from a stationary wave from $U_1\in\CC_+$.
 Observe that both $U_1^o$ and $U_2^{\#}$ belong to $\WW_3(U_1)$. Whenever $U_1^o$ lies above  $U^{\#}$ on $\WW_3(U_1)$, there are three distinct solutions. Otherwise, there is at most one solution.  See Figure \ref{Riem2-multiple}.

\begin{theorem}[Riemann problem for the shallow water equations]
 \label{theo33} Given a left-hand state $U_L\in G_2$.
\begin{itemize}

\item[(a)] {\bf Existence.} The Riemann problem \rf{1.1}--\rf{1.3} admits a solution  if $K^o$ lies below or on the curve $\WW_2^B(U_R)$, and that if $\WW_2^B(U_R)$ intersects with $\WW_1^a(U_L)$ at some point $U_M^o\in G_2^-$, then   $\bar\ld_2(U_M^o,U_R)\ge 0$.

\item[(b)] {\bf Regime of uniqueness.} The Riemann problem \rf{1.1}--\rf{1.3} has at most one solution if

-- either $a_L\le a_R$;

-- or $a_L>a_R$, $h_1^o\ge h_2^{\#}$, where $U_2$ is defined in \rf{3.8};

-- or $a_L>a_R$, $h_1^o< h_2^{\#}$, and the states $U_1^o$ and  $U_2^{\#}$ are located on the same side with respect to the curve $\WW_2^B(U_R)$.

\item[(c)] {\bf Multiple  solutions.} If $a_L>a_R$, $h_1^o< h_2^{\#}$, and $U_1^o$ lies above the curve $\WW_2^B(U_R)$  and $U_2^{\#}$ lies below the curve $\WW_2^B(U_R)$,  then the Riemann problem \rf{1.1}--\rf{1.3} has three solutions.
\end{itemize}
\end{theorem}

\noindent{\it Example.} We provide  some numerical experiments computing $h_1^o, h_2^{\#}$ to illustrate the cases of Theorem \ref{theo33} (see Tables B1, B2, and B3). All experiments give the same result: $h_1^o>h_2^{\#}$.

$$
\framebox{
\begin{tabular}{l l c l c l c}
{\bf Table B1} \\
\hline
   States          &\vline  &  $U_L$   &\vline &   $U_1^o$ & \vline &   $U_2^{\#}$ \\
\hline
Water Height $h$   & \vline &     3  & \vline &      1.819500899801235  &\vline &   1.768961248574716  \\
Velocity  $u$      & \vline &     0.5    & \vline &     3.032474262659020  &\vline &      3.119112786658156  \\
Bottom Level  $a$  & \vline &     1.1    & \vline &      1.000000000000000  &\vline &      1.000000000000000  \\
\end{tabular}
}
\label{TableB1}
$$
$$
\framebox{
\begin{tabular}{l l c l c l c}
{\bf Table B2} \\
\hline
   States          &\vline  &  $U_L$   &\vline &   $U_1^o$ & \vline &   $U_2^{\#}$ \\
\hline
Water Height $h$   & \vline &     3  & \vline &          1.707571536932233 &\vline &     1.656818524474798  \\
Velocity  $u$      & \vline &     0.1    & \vline &         2.901359698616083  &\vline & 2.990236508448978\\
Bottom Level  $a$  & \vline &     1.1    & \vline &         1.000000000000000  &\vline &         1.000000000000000 \\
\end{tabular}
}
\label{TableB2}
$$
$$
\framebox{
\begin{tabular}{l l c l c l c}
{\bf Table B3} \\
\hline
   States          &\vline  &  $U_L$   &\vline &   $U_1^o$ & \vline &   $U_2^{\#}$ \\
\hline
Water Height $h$   & \vline &     3  & \vline &         3.187878980786353 &\vline &      2.574902018055705 \\
Velocity  $u$      & \vline &     1    & \vline &        1.969891931767155  &\vline & 2.438841182952260\\
Bottom Level  $a$  & \vline &     2    & \vline &         1.000000000000000  &\vline &         1.000000000000000 \\
\end{tabular}
}
\label{TableB3}
$$

\begin{rmk}
We conjecture that $h_1^o>h_2^{\#}$.
If this conjecture holds, then Theorem \ref{theo33} implies that the Riemann problem always has at most one solution for $U_L\in G_2$.
\end{rmk}

\subsection{Remark on the continuous dependence of solutions}

As seen in the previous subsections, the construction of Riemann solutions is
 based on a given left-hand state $U_L$. The Riemann problem for \rf{1.1}--\rf{1.2} admits up to three solutions for data in certain regions, which implies that the initial-value problem for \rf{1.1}--\rf{1.2} is ill-posed.  However, connectivity between the types of Riemann solutions helps to determine the continuous dependence of the set of solutions on Riemann data.  
Observe that for each construction (A) and (B), the general 
structure of solutions changes continuously when $U_R$ changes and one evolves from one case to another. For example, Construction (A1) changes continuously to (A2), and 
(A2) itself changes continuously to (A3). Similar remarks hold for the cases (B1), (B2), and (B3), as observed earlier 
about the continuity of the wave patterns. Thus, the set of solutions ``globally''  
depends continuously on the right-hand side $U_R$ for each case $U_L\in G_1\cup \CC_+$ and $U_L\in G_2$. In order to show that the set of solutions depends continuously on Riemann data, we need only to check
that when $U_L$ moves from $G_1$ to $G_2$, the change in the structure of solutions is continuous as well. But this fact also holds true, since when $U_L$ tends to $\CC_+$ on each side, the solutions in Constructions (A1) and (B1) approach each other,  and the solution of Constructions (A3) and (B3) approach
 each other as long as the solutions make sense. If $a_L\ge a_R$, these solutions eventually coincide on $\CC_+$.

\section{A Godunov-type algorithm}

\subsection{A well-balanced quasi-conservative scheme}

Given a uniform time step $\Delta t$, and a spatial mesh size
$\Delta x$, setting $x_i=i\Delta x, i\in\ZZ$, and $ t^n = n\Delta
t, n\in \NN$, we denote $U_{i}^n$ to be an approximation of the exact value $U(x_i,t^n)$.
Set
$$
U=\left(\begin{array}{c} h\\
hu\\
\end{array}\right),\quad F(U)=\left(\begin{array}{c} hu\\
h(u^2 + g\frac{h}{2})\\
\end{array}\right),\quad S(U) = \left(\begin{array}{c} 0\\
- gh\\
\end{array}\right)\pt_x a.
$$
The system \rf{1.1}--\rf{1.2} can be written in the compact form
 \begin{equation}
\pt_tU  + \pt_x F(U)=S(U)\pt_x a, \quad t>0,x\in\RR.
\label{4.1}
\end{equation}
Let us be given the initial condition
\begin{equation}
U(x,0)=U_0(x),\quad x\in\RR,
\label{4.2}
\end{equation}
Then, the discrete initial values  $U_i^0$ are given by
\begin{equation}
U_i^{0}= {1\over\Delta x}\int_{x_{i-1/2}}^{x_{i+1/2}}U_0(x) dx.
\label{4.3}
\end{equation}
Suppose  $U^n$ is known and $U^n$  is constant on each interval $(x_{i-1/2},x_{i+1/2})$ v�i $i\in \ZZ$. On each cell $(x_{i-1},x_i)$ we determine the exact solution to the Riemann problem for
\begin{equation}
\pt_tU(x,t)  + \pt_x F(U(x,t))=S(U)\pt_x a,\quad \textrm{on } \RR\times (t^n,t^{n+1}],
\label{4.4}
\end{equation}
subject to the initial condition
\begin{equation}
U(x,t^n) = \left\{\begin{array}{ll} U_{i-1}^n,&\quad  x<x_{i-1/2},\\
U_{i}^n,& x>x_{i-1/2}.\\
\end{array}\right.
\label{4.5}
\end{equation}
Denote this solution by  $U(x,t;U_{i-1}^n, U_{i}^n)$.
Use these solutions of the local Riemann problems we define the function $V$ by
$$
V(x,t):=\left\{\begin{array}{ll}U(x,t;U_{i-1}^n, U_{i}^n),&\quad x_{i-1/2}< x\le x_{i}, t^n\le t\le t^{n+1},\\
U(x,t;U_{i}^n, U_{i+1}^n),&\quad x_{i}< x\le x_{i+1/2}, t^n\le t\le t^{n+1}.\\
\end{array}\right.
$$
As for the initial values, we have to ensure that the approximation $U_i^{n+1}$ at time $t^{n+1}$ is constant on $(x_{i-1/2},x_{i+1/2})$ for all $i\in \ZZ$.
Therefore, we define the new value  $U_i^{n+1}$ at the time $t=t^{n+1}$ by
\begin{equation}
U_i^{n+1}= {1\over\Delta x}\int_{x_{i-1/2}}^{x_{i+1/2}}V(x,t^{n+1})dx.
\label{4.6}
\end{equation}
This means $U_i^{n+1}$ is the mean value of $V$ on $(x_{i-1/2},x_{i+1/2})$ and thus contains parts of
$U(x,t;U_{i-1}^n, U_{i}^n)$ and $U(x,t;U_{i}^n, U_{i+1}^n)$. To ensure that the solutions of two consecutive local Riemann problems do not coincide, we assume that the following CFL (Courant, Friedrichs, Lewy)  condition holds:
$$
{\Delta t\over \Delta x}\max|\ld_i|  \le {1\over 2},
$$
where $\ld_i$ denote the eigenvalues of $DF(U)$.

Suppose now $V$ is an exact solution on $(x_{i-1/2},x_{i+1/2})$. Since the $a$-component is constant in $(x_{i-1/2},x_{i+1/2})$, the right-hand side of \rf{1.1} vanishes for $V$. Thus, the standard Godunov scheme is in {\it quasi-conservative form}:
\begin{equation}
U_i^{n+1}= U_i^{n}-{\Delta t\over\Delta x} (F(U(x_{i+1/2}-,t^{n+1};U_{i}^n, U_{i+1}^n))-F(U(x_{i-1/2}+,t^{n+1};U_{i-1}^n, U_{i}^n))).
\label{4.7}
\end{equation}
One might think that in the scheme \rf{4.7} the source term is incorporated into the local Riemann problem.

The Godunov scheme \rf{4.7} is capable of capturing exactly equilibria. Therefore \rf{4.7} is a {\it well-balanced} scheme. In fact, let us be given the initial data $U^0$ to be equilibrium states of a stationary wave. Then, on each cell $x_{i-1/2}<x<x_{i+1/2}, t^n<t\le t^{n+1}$ the exact Riemann solution is constant. Thus, $U(x_{i+1/2}-,t^n;U_{i}^n, U_{i+1}^n)=U(x_{i-1/2}+,t^n;U_{i-1}^n, U_{i}^n)$ and so $U_i^{n+1}= U_i^{n}=U^0_i$ for all $i\in \ZZ$ and $n\ge 0$.
When there are multiple Riemann solutions, any of them can be selected and we still  
obtain a deterministic scheme, according to Theorem \ref{theo32}.

\subsection{Numerical  Riemann solver}

Given any Riemann data $(U_L,U_R)$, denote by $U(x,t;U_L,U_R)$ the Riemann solution corresponding to the Riemann data $(U_L,U_R)$. To build the Godunov scheme \rf{4.7} we will specify the values $U(0\pm,\Delta t;U_L,U_R)$ for an arbitrary and fixed  number $\Delta t>0$.

\noindent\underline{\it Riemann solver (A1).} We present a computing strategy  for Riemann solutions \rf{3.1} as follows.
\begin{itemize}
\item[(i)] The state $U_L^o=(h_L^o,u_L^o,a_R)$: $h_L^o=\bar h(a_R)=h_1(a_R)$ , where $h_1(a_R)$ is the smaller root of the nonlinear equation \rf{2.13}, described by Lemma \ref{lem31}, and can be computed using Lemma \ref{lem23}. $u_L^o=u_Lh_L/h_L^o$.
\item[(ii)] The state $U_M=(h_M,u_M,a_R)$ is the intersection point of the wave curves $\WW_1(U_L^o)$ and $\WW_2^B(U_R)$, see \rf{2.10}. Equating the $u$-component for these two curves leads to a strictly increasing and strictly convex function in $h$. Thus, the $h$-component of the intersection point $h_2$ can be computed using the Newton's method.
\end{itemize}
The Riemann solver (A1) relying on Construction (A1) yields
\begin{equation}
\begin{array}{ll}
&U(0-,\Delta t;U_L,U_R)=U_L,\\
&U(0+,\Delta t;U_L,U_R)=U_L^o.\\
\end{array}
\label{4.8}
\end{equation}
This implies that the Godunov scheme \rf{4.7} using the Riemann solver 1 becomes
\begin{equation}
U_i^{n+1}= U_i^{n}-{\Delta t\over\Delta x} (F(U_{i}^n) - F((U_{i-1}^n)^o)),
\label{4.9}
\end{equation}
where $U^o$ defined as in \rf{4.8}.

\noindent\underline{\it Riemann solver  (A2).}

The states of the Riemann solution \rf{3.3} can be found as follows.
\begin{itemize}
\item[(1)] The state $U_M=(h_M,u_M,a_R)$  is determined by
$$\{U_M\}=\WW_3(U_L)\cap \WW_2^B(U_R).
$$
\item[(2)] The states $U_1=(h_1,u_1,a_1), U_2=(h_2,u_2,a_1)$ are determined by using the corresponding
``half-way" shifting in $a$ component from the stationary contact from $U_L$ to $U_1$ and the  stationary contact from $U_2$ to $U_M$, and using the fact that $U_2=U_1^{\#}$, (see Lemma \ref{lem31}):
    \begin{equation}
\begin{array}{ll}
&a_1=a_L+{u_L^2-u_1^2\over 2g}+h_L-h_1=a_M+{u_M^2-u_2^2\over 2g}+h_M-h_2,\\
& u_1=\frac{u_Lh_L}{h_1},\\
&h_2=\frac{-h_1+\sqrt{h_1^2+8h_1u_1^2/g}}{2},\\
& u_2=\frac{u_Mh_M}{h_2}.\\
\end{array}
\label{4.10}
\end{equation}
\end{itemize}
It is not difficult to check that the system \rf{4.10} can yield a scalar equation for $h_1$.
The Riemann solver (A2) relying Construction (A2) gives
\begin{equation}
\begin{array}{ll}
&U(0-,\Delta t;U_L,U_R)=U_L,\\
&U(0+,\Delta t;U_L,U_R)=U_M=U_{M}(U_L,U_R).\\
\end{array}
\label{4.11}
\end{equation}
This implies that the Godunov scheme \rf{4.7} using the Riemann solver 2 becomes
\begin{equation}
U_i^{n+1}= U_i^{n}-{\Delta t\over\Delta x} (F(U_{i}^n) - F(U_M(U_{i-1}^n,U_i^n))),
\label{4.12}
\end{equation}
where $U_M(U_{i-1}^n,U_i^n))$ is defined as in \rf{4.11}, i.e.
$$
\{U_M(U_{i-1}^n,U_i^n))\}= \WW_3(U_{i-1}^n)\cap \WW_2^B(U_{i}^n).
$$
Since $U_M$ plays a key role in this Riemann solver, we sketch a computing algorithm for $U_M$ as follows. First we observe that if $U_R$ lies below the curve $\WW_3(U_L)$ in the $(h,u)$-plane, then $U_M$ is the intersection point of $\WW_3(U_L)$ and $\Scal_2^B(U_R)$. Otherwise, $U_M$ is the intersection point of $\WW_3(U_L)$ and $\Rcal_2^B(U_R)$.
Thus, we find:
\begin{itemize}
\item[(i)] (Arrival by a $2$-shock) If $h_Ru_R-h_Lu_L<0$ then $h_M$ is the root of the equation
\begin{equation}
G_1(h):={h_Lu_L\over h}-\Big(u_R+(h-h_R)\sqrt{{g\over 2}\big({1\over h}+{1\over h_R}\big)}\Big)=0.
\label{4.13}
\end{equation}
\item[(ii)] (Arrival by a $2$-rarefaction wave) Otherwise, $h_M$ is the root of the equation
\begin{equation}
G_2(h):={h_Lu_L\over h}-\big(u_R+2\sqrt{g}(\sqrt{h}-\sqrt{h_R})\big)=0.
\label{4.14}
\end{equation}
\end{itemize}
It is easy to see that both functions $G_1, G_2$ defined by \rf{4.13} and \rf{4.14} are strictly convex. Moreover,
we have
$$
\begin{array}{ll}G_1'(h)&=-{h_Lu_L\over h^2}-\sqrt{g\over 2}\Big({1\over h}+{2\over h_R}+{h_R\over h^2}\Big)\Big(\frac{1}{2\sqrt{1/h+1/h_R}}\Big)<0,\\%\frac{1}{h}+\frac{1}{h_R}
G_2'(h)&=-{h_Lu_L\over h^2}-\sqrt{g\over h}<0,
\end{array}$$
for all $h>0$. Thus, the Newton method can be applied for both equations \rf{4.13} and \rf{4.14} with any starting point.

\noindent\underline{\it Riemann solver (A3).}

Let us consider Construction (A3) and let
$$
A=U_L^{\#},\quad \{B\}=\WW_1(U_L)\cap \WW_2^B(U_R).
$$
It is easy to see that $U_M=(h_M,u_M,a_L)$ lies on $\WW_1(U_L)$ between $A$ and $B$.
We propose a procedure similar to the Bisection method to compute the states of the elementary waves of the Riemann solution \rf{3.5} as follows.
 We use the equation of $\WW_2^B(U_R): \Phi_2(U;U_R)=0$, defined by \rf{2.10}, as a test condition: for $U$ above $\WW_2^B(U_R)$, $\Phi_2(U;U_R)>0$ and for $U$ below it, $\Phi_2(U;U_R)<0$.
 Using a stationary jump from any state $U$ on the wave pattern of $\WW_1(U_L)$ between $A$ and $B$ to a state $U^o$ shifting $a$ from $a_L$ to $a_R$. Then, we have
$$
\Phi_2(A;U_R)\cdot \Phi_2(B;U_R)<0.
$$

\underline{\bf Algorithm 1:}

\underline{Step 1:} An estimate for $h_M$ is given by
$$
h_M={h_A+h_B\over 2},
$$
$U_M=(h_M,u_M,a_L)\in \WW_1(U_L)$, so $u_M$ is computed using the equation \rf{2.6} with $U_0=U_L$.

\underline{Step 2:}

\begin{itemize}
\item[(a)]  If $\Phi_2(A;U_R) \cdot \Phi_2(U_M;U_R)<0$, then set $B=U_M$ and return to Step 1;

\item[(b)]  If $\Phi_2(A;U_R) \cdot \Phi_2(U_M;U_R)>0$, then set $A=U_M$ and return to Step 1;

\item[(c)]  If $\Phi_2(A;U_R) \cdot \Phi_2(U_M;U_R)=0$, terminate the computation.
\end{itemize}

We can still use an alternative algorithm using the value of $a$-component as a convergence condition, as follows.

\underline{\bf Algorithm 2:}

\underline{Step 1:} Let $A=U_L^{\#}$ and $B$ is the intersection point of $\WW_1(U_L)$ and $\{u=0\}$. An estimate for $h_M$ is given by
$$
h_M={h_A+h_B\over 2},
$$
and $u_M$ is estimated using the equation \rf{2.6}, so an estimate of $U_M$ is $U_M=(h_M,u_M,a_L)\in \WW_1(U_L)$.
An estimate for $U_M^o$ is given by
$$
\{U_M^o\}=\WW_3(U_M)\cap \WW_2^B(U_R).
$$
Determine the change in $a$-component for the stationary wave between $U_M$ and $U_M^o$ (see \rf{2.12})
$$
a=a_L+{u_M^2-(u_M^o)^2\over 2g}+h_M-h_M^o.
$$
\underline{Step 2:}

\begin{itemize}
\item[(a)]  If $a-a_R<0$, then set $h_A=h_M$ and return to Step 1;

\item[(b)]  If $a-a_R>0$, then set $h_B=h_M$ and return to Step 1;

\item[(c)]  If $a-a_R=0$, stop the computation.
\end{itemize}

The Riemann solver (A3) relying on Construction (A3) yields
\begin{equation}
\begin{array}{ll}
&U(0-,\Delta t;U_L,U_R)=U_M=U_M(U_L,U_R),\\
&U(0+,\Delta t;U_L,U_R)=U_M^o=(U_M(U_L,U_R))^o.\\
\end{array}
\label{4.15}
\end{equation}
This implies that the Godunov scheme \rf{4.7} using the Riemann solver 3 becomes
\begin{equation}
U_i^{n+1}= U_i^{n}-{\Delta t\over\Delta x} (F(U_M(U_i^n,U_{i+1}^n)) - F(U_M^o(U_{i-1}^n,U_i^n))),
\label{4.16}
\end{equation}
where $U_M(U_i^n,U_{i+1}^n)$ and $U_M^o(U_{i-1}^n,U_i^n))$ are defined as in \rf{4.15}.

\noindent\underline{\it Riemann solver (B1).}
The Riemann solver (B1) relying on  Construction (B1) gives
\begin{equation}
\begin{array}{ll}
&U(0-,\Delta t;U_L,U_R)=U_1=\Big(\big({u_L\over 3\sqrt{g}}+{2\over 3}\sqrt{h_L}\big)^2, {1\over 3}u_L+{2\over 3}\sqrt{gh_L},a_L\Big):=U_{L,+}\in \CC_+,\\
&U(0+,\Delta t;U_L,U_R)=U_2:=U_{L,+o}\in G_1,\\
\end{array}
\label{4.17}
\end{equation}
If $a_L\ge a_R$, then
$$
\begin{array}{ll}
&U_1=\Big(\big({u_L\over 3\sqrt{g}}+{2\over 3}\sqrt{h_L}\big)^2, {1\over 3}u_L+{2\over 3}\sqrt{gh_L},a_L\Big):=U_{L,+}\in \CC_+,\\
&U_2:=U_{L,+o}\in G_1,\\
\end{array}
$$
where $U_{L,+o}\in G_1$ is the state resulted by a stationary contact from $U_{L,+}\in\CC_+$.
This implies that the Godunov scheme \rf{4.7} using the Riemann solver (B1) becomes
\begin{equation}
U_i^{n+1}= U_i^{n}-{\Delta t\over\Delta x} (F(U_{i,+}^n) - F(U_{i-1,+o}^n)).
\label{4.18}
\end{equation}
If  $a_L<a_R$, then $U_2\in\CC_+$. The computing of $U_1$ and $U_2$ can be done similarly as in the Riemann solver (A3).

\noindent\underline{\it Riemann solver (B2).}
The Riemann solver (B2) relying Construction (B2) gives
\begin{equation}
\begin{array}{ll}
&U(0-,\Delta t;U_L,U_R)=U_1=\Big(\big({u_L\over 3\sqrt{g}}+{2\over 3}\sqrt{h_L}\big)^2, {1\over 3}u_L+{2\over 3}\sqrt{gh_L},a_L\Big):=U_{L,+}\in \CC_+,\\
&U(0+,\Delta t;U_L,U_R)=P=P(U_L,U_R)\in \WW_2^B(U_R)\cap \WW_3(U_{L,+}).\\
\end{array}
\label{4.19}
\end{equation}
This implies that the Godunov scheme \rf{4.7} using the Riemann solver 2 becomes
\begin{equation}
U_i^{n+1}= U_i^{n}-{\Delta t\over\Delta x} (F(U_{i,+}^n) - F(P(U_{i-1}^n,U_i^n))),
\label{4.20}
\end{equation}
where
$$
P(U_{i-1}^n,U_i^n))=\WW_2^B(U_i^n)\cap \WW_3(U_{i-1,+}^n).
$$

\noindent\underline{\it Riemann solver (B3).}
The Riemann solver (B3) relying on Construction (B3) yields
\begin{equation}
\begin{array}{ll}
&U(0-,\Delta t;U_L,U_R)=U_M=U_M(U_L,U_R),\\
&U(0+,\Delta t;U_L,U_R)=U_M^o=(U_M(U_L,U_R))^o.\\
\end{array}
\label{4.21}
\end{equation}
This implies that the Godunov scheme \rf{4.7} using the Riemann solver 3 becomes
\begin{equation}
U_i^{n+1}= U_i^{n}-{\Delta t\over\Delta x} (F(U_M(U_i^n,U_{i+1}^n)) - F(U_M^o(U_{i-1}^n,U_i^n))),
\label{4.22}
\end{equation}
where $U_M(U_i^n,U_{i+1}^n)$ and $U_M^o(U_{i-1}^n,U_i^n))$ are defined as in \rf{4.21}. It is easy to see that the determinations of the states $U(0-,\Delta t;U_L,U_R)$ and $U(0+,\Delta t;U_L,U_R)$ by Solvers (A3) and (B3) are the same.

The computing strategies for Solvers (B1), (B2), and (B3) are similar to those of Solvers (A1), (A2), and (A3).

\subsection{The Godunov algorithm}

It is natural to ask which Riemann solvers should be taken in the Godunov scheme.
As seen earlier, in the regions where there are possibly multiple solutions one can select {\sl any} Riemann solution. 
Three extreme cases can be distinguished by preferring one particular Riemann solver whenever it is available. For example, we can decide to always select solutions with stationary contact wave in the same region as the left-hand state, that is, the solvers (A1) and (B3). This selection leads us to a deterministic
 algorithm for designing a corresponding Godunov scheme, as now  described.

\bigskip

\noindent{\it Building Godunov Scheme Algorithm preferring solvers (A1) and (B3).} Let $U_L=U_i^n$ and $U_R=U_{i+1}^n$.

\bigskip

If $\ld_1(U_L)\le 0$

 \quad   If $\Phi_2(U_L^{o\#};U_R)<0$

          \qquad  Use Solver (A1)

  \quad   elseif $\Phi_2(U_L^{\#o};U_R)<0$

           \qquad Use Solver (A2)

   \quad  else

           \qquad Use Solver (A3)

   \quad   end

else

    \quad If $\Phi_2(U_1^o;U_R)>0$

           \qquad Use Solver (B3)

  \quad   elseif $\Phi_2(U_2^{\#};U_R)>0$

         \qquad   Use Solver (B2)

  \quad   else

     \qquad    Use Solver (B1)

   \quad   end

end

\section{Numerical experiments (I). The non-resonant regime}

We now numerically investigate our Riemann solver and Godunov method and present several numerical tests.
For each test we consider
 the errors between the exact Riemann solution and the approximate solution by the Godunov scheme \rf{4.7} for $x\in [-1, 1]$ with different mesh sizes corresponding to
 $500, 1000, 2000$ points. In this section as well as in the next section, we plot the solution at
  the time $t=0.1$, and use the  stability condition
$$
 CFL=0.75. 
 $$
The algorithm for selecting the Riemann solvers is the one described at the end of the last section, unless 
indicated otherwise. 

\subsection{Test 1}

\begin{figure}[h]
\begin{center}
  \includegraphics[width=10cm]{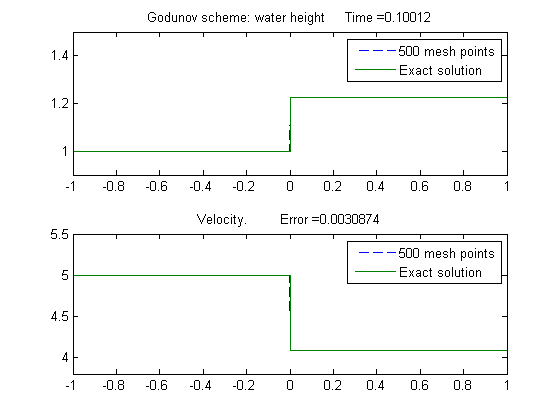}\\
  \caption{Test 1: A stationary contact wave is captured exactly by Godunov method
using our exact Riemann solver.}\label{T1}
\end{center}
\end{figure}

 This test indicates that the Godunov method is capable of maintaining equilibrium states. Let
\begin{equation}
U_0(x)= \left\{\begin{array}{ll}
&U_L=(h_L,u_L,a_L),\quad x<0\\
&U_R=(h_R,u_R,a_R),\quad x>0,\\
\end{array}\right.\label{5.1}
\end{equation}
where $U_L=(1,5,1)$ and $U_R=(1.223655890827479, 4.086116070277590,1.2)$.
It is not difficult to check that the Riemann problem with initial data \rf{5.1} admits  a stationary contact between these equilibrium states:
\begin{equation}
 U(x,t)=U_0(x),\quad x\in\RR, t>0.
\label{5.2}
\end{equation}
Figure \ref{T1} shows that the stationary contact is well captured by Godunov method using our exact Riemann solver for $x\in [-1, 1]$ with $500$ mesh points and at time $t=0.1$.

\subsection{Test 2}

\begin{figure}
  \includegraphics[width=10cm]{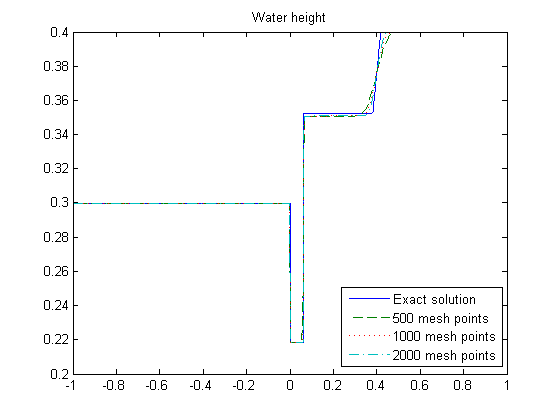}\\
  \caption{Test 2.  Water height with different mesh sizes.}\label{T2h}
\end{figure}

\begin{figure}
  \includegraphics[width=10cm]{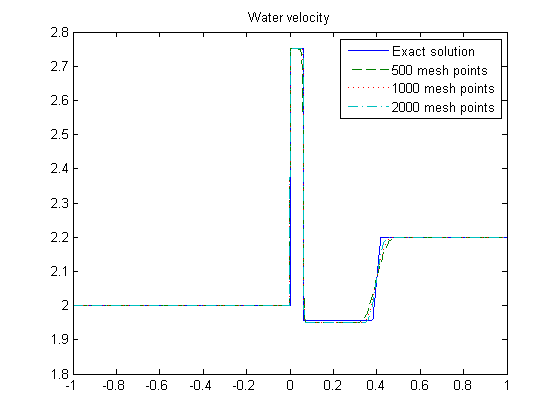}\\
  \caption{Test 2.  Water velocity with different mesh sizes.}\label{T2u}
\end{figure}

We now approximate a non-stationary Riemann solution with data $U_L,U_R\in G_1$.  Precisely, we consider the Riemann problem \rf{1.1}--\rf{1.3} with data
\begin{equation}
U_0(x)= \left\{\begin{array}{ll}
&U_L=(h_L,u_L,a_L)=(0.3,2,1.1)\in G_1,\quad x<0,\\
&U_R=(h_R,u_R,a_R)=(0.4,2.2,1)\in G_1,\quad x>0.\\
\end{array}\right.\label{5.3}
\end{equation}
The Riemann problem \rf{1.1}--\rf{1.2} with the initial data \rf{5.3} admits the solution described by Construction 1, where
$$
\begin{array}{ll}&U_1=(0.21815897, 2.750288, 1), 
\qquad \quad 
U_2=(0.35252714, 1.9572394,1).\\
\end{array}$$
The errors for Test 2 are reported in the following table
%Table \ref{5.4}.
$$
%\begin{equation}
\framebox{
\begin{tabular}{c l c l c}
% column 1 entry & column 2 entry ... & column n entry \\
$N$    &\vline &   $||U_h^\textrm{\scriptsize C} - U||_{L^1}$ &\\
\hline
500   & \vline &   0.012644   \\
1000   & \vline &    0.0087928  \\
2000   & \vline &   0.0063773  \\
\end{tabular}
} 
%\label{5.4}
$$
%\end{equation}
and Figures \ref{T2h} and \ref{T2u} and the table above
%Table \rf{5.4} 
show that approximate solutions are closer to the exact solution when the mesh size gets smaller.

\subsection{Test 3}

\begin{figure}
  \includegraphics[width=10cm]{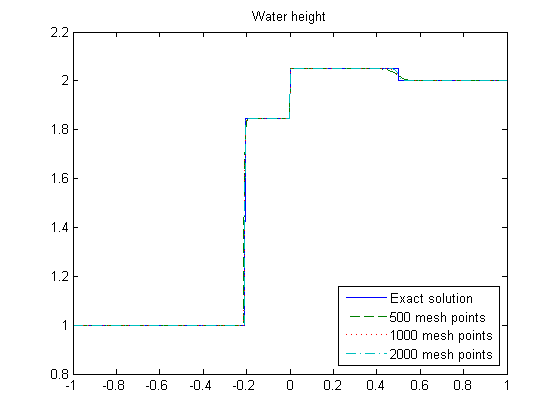}\\
  \caption{Test 3.  Water height with different mesh sizes.}\label{T3h}
\end{figure}

\begin{figure}
  \includegraphics[width=10cm]{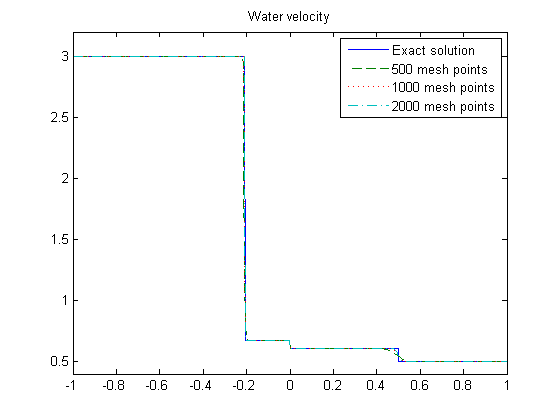}\\
  \caption{Test 3. Water velocity with different mesh sizes.}\label{T3u}
\end{figure}

In this test, we will approximate a non-stationary Riemann solution with Riemann data $U_L,U_R\in G_2$.  Precisely, we consider the Riemann problem \rf{1.1}--\rf{1.3} with data
\begin{equation}
U_0(x)= \left\{\begin{array}{ll}
&U_L=(h_L,u_L,a_L)=(1,3,1.2)\in G_2,\quad x<0,\\
&U_R=(h_R,u_R,a_R)=(2,0.5,1)\in G_2,\quad x>0.\\
\end{array}\right.\label{5.5}
\end{equation}
This problem admits the solution described by Construction (A3), where
$$
\begin{array}{ll}&U_M=(1.8452179,0.67672469,1.2), 
\qquad \quad 
U_M^o=(2.0496463,0.60922927,1).\\
\end{array}
$$
Precisely, the solution starts with a 1-shock from $U_L$ to $U_M$, followed by a stationary contact from $U_M$ to $U_M^o$, and then arrived at $U_R$ from $U_M^o$ by a 2-shock wave. 
The errors for Test 3 are reported in the following table: 
%Table \ref{5.6}.
\begin{equation}
\framebox{
\begin{tabular}{c l c l c}
% column 1 entry & column 2 entry ... & column n entry \\
$N$    &\vline &   $||U_h^\textrm{\scriptsize C} - U||_{L^1}$ &\\
\hline
500   & \vline &  0.01813  \\
1000   & \vline &   0.0076434 \\
2000   & \vline &    0.0035277\\
\end{tabular}
} \label{5.6}
\end{equation}
Figures \ref{T3h} and \ref{T3u} and the table above
%Table \ref{5.6} 
show that approximate solutions are closer to the exact solution when the mesh size gets smaller.
All of our tests presented so far exhibit a convergence of approximate solutions to the exact solution.

\section{Numerical experiments (II). The resonance regime}

In the following, we will consider the cases where the Riemann data on the different sides with respect to $\CC_+$:
\begin{itemize}
\item[(i)]  $U_L\in G_1$ and $U_R\in G_2$;
\item[(ii)]  $U_L\in G_2$ and $U_R\in G_1$.
\end{itemize}
 The solution is evaluated for $x\in [-1, 1]$ with $500$ points and at time $t=0.1$. We take also
 $$
 CFL=0.75.
 $$

\subsection{Test 4}

\begin{figure}
  \includegraphics[width=10cm]{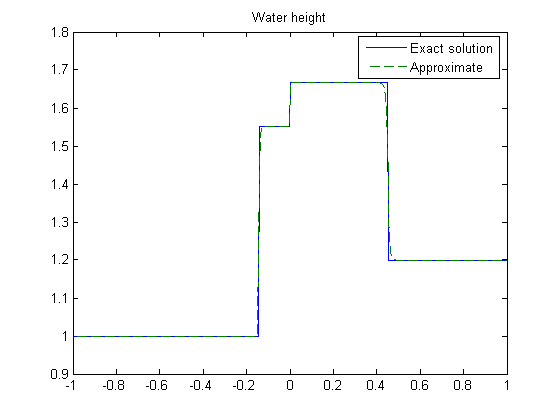}\\
  \caption{Test 4 (Resonant case).
 Water height of the solution  \rf{3.5}.}\label{A3-Resonance1h}
\end{figure}

\begin{figure}
  \includegraphics[width=10cm]{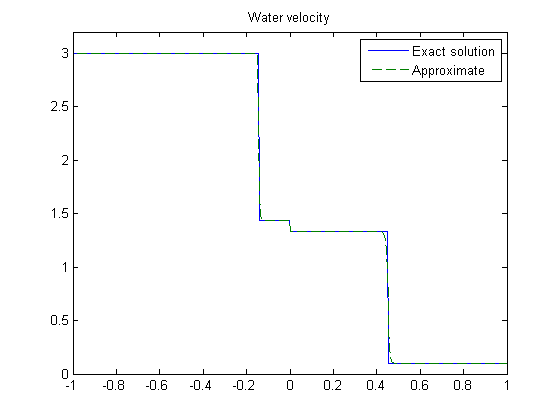}\\
  \caption{Test 4 (Resonant case).
 Water velocity of the solution \rf{3.5}.}\label{A3-Resonance1u}
\end{figure}

In this test, $a_L>a_R$, and there is a unique solution.  We consider the Riemann problem
 \rf{1.1}--\rf{1.3} with data
\begin{equation}
U_0(x)= \left\{\begin{array}{ll}
&U_L=(h_L,u_L,a_L)=(1, 3, 1.1)\in G_1,\quad x<0,\\
&U_R=(h_R,u_R,a_R)=(1.2,0.1,1)\in G_2,\quad x>0.\\
\end{array}\right.\label{6.1}
\end{equation}
We have
$$
h_L^{o\#}=1.042865405801653 <  h_L^{\#o} = 1.213385283426733.
$$
Thus,  the  problem \rf{1.1}, \rf{1.2}, \rf{6.1} admits a unique solution of the form \rf{3.1},  according to Theorem \rf{theo32}, where
$$
%\begin{equation}
U_M=(1.5521168,1.4328264,1.1),\quad U_M^o=(1.665941, 1.3349296, 1).
%\label{6.2}
$$
%s\end{equation}
Figures \ref{A3-Resonance1h}-\ref{A3-Resonance1u}  show that  the Godunov scheme gives  good  approximate solutions to the exact solution in this resonance case.

 \subsection{Test 5}

\begin{figure}
  \includegraphics[width=10cm]{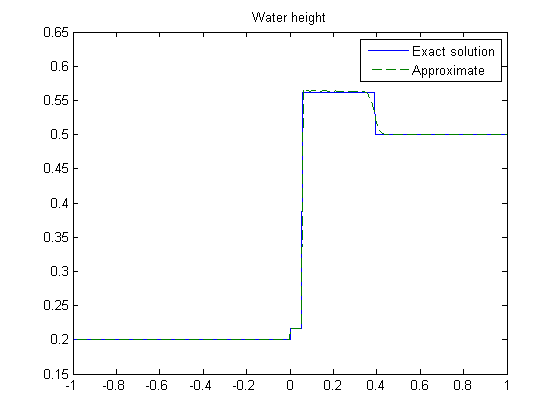}\\
  \caption{Test 5 (Resonant case).
 Water height of the solution  \rf{3.1}.}\label{A1c-Resonance1h}
\end{figure}

\begin{figure}
  \includegraphics[width=10cm]{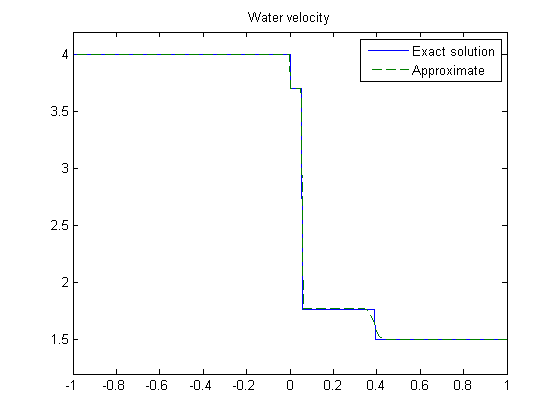}\\
  \caption{Test 5  (Resonant case).
 Water velocity of the solution \rf{3.1}.}\label{A1c-Resonance1u}
\end{figure}

In this test, $a_L<a_R$, and there is a unique solution.  We consider the Riemann problem
 \rf{1.1}--\rf{1.3} with data
\begin{equation}
U_0(x)= \left\{\begin{array}{ll}
&U_L=(h_L,u_L,a_L)=(0.2, 4, 1)\in G_1,\quad x<0,\\
&U_R=(h_R,u_R,a_R)=(0.5,1.5,1.1)\in G_2,\quad x>0,\\
\end{array}\right.\label{6.3}
\end{equation}
and
$$
\begin{array}{ll}
&U_L^{o\#}=(0.677264819960833,1.181221844722815),\\
&\Phi_2(U_L^{o\#};U_R)=  -1.050411375011095<0,\\
&U_L^{\#o}=(0.581828763814630, 1.374974992221044),\\
&\Phi_2(U_L^{\#o};U_R)=-0.474326705580410<0.\\
\end{array}
$$
This Riemann problem admits a unique solution of the form \rf{3.1},  according to Theorem \ref{theo32}, where
$$
%\begin{equation}
U_L^o=(0.21591647, 3.7051366,1.1),\quad U_M=(0.56185289, 1.7661913, 1.1).
%\label{6.4}
$$
%\end{equation}
Figures \ref{A1c-Resonance1h} to \ref{A1c-Resonance1u}  show that  the Godunov scheme gives  good  approximate solutions to the exact solution in this resonance case.

 \subsection{Test 6}

Next, we provide a test for the case there are multiple solutions. So we consider the Riemann problem 
 \rf{1.1}--\rf{1.3} with data
\begin{equation}
U_0(x)= \left\{\begin{array}{ll}
&U_L=(h_L,u_L,a_L)=(0.2, 5, 1)\in G_1,\quad x<0,\\
&U_R=(h_R,u_R,a_R)=(0.75904946,1.3410741, 1.2)\in G_2,\quad x>0.\\
\end{array}\right.\label{6.5}
\end{equation}

\begin{figure}
  \includegraphics[width=10cm]{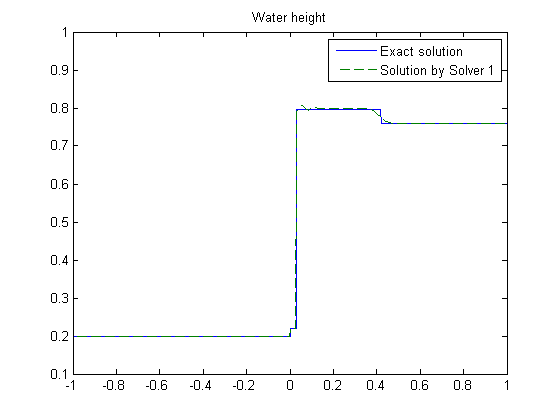}\\
  \caption{Test 6 (Resonant case - multiple solutions).
 Water height of the solution \rf{3.1} - preferred Solver (A1).}\label{Resonance1h}
\end{figure}

\begin{figure}
  \includegraphics[width=10cm]{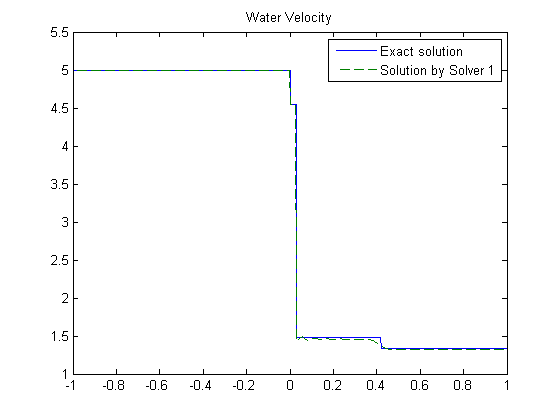}\\
  \caption{Test 6 (Resonant case - multiple solutions).
  Water velocity of the solution \rf{3.1} - preferred Solver (A1).}\label{Resonance1u}
\end{figure}

\begin{figure}
  \includegraphics[width=10cm]{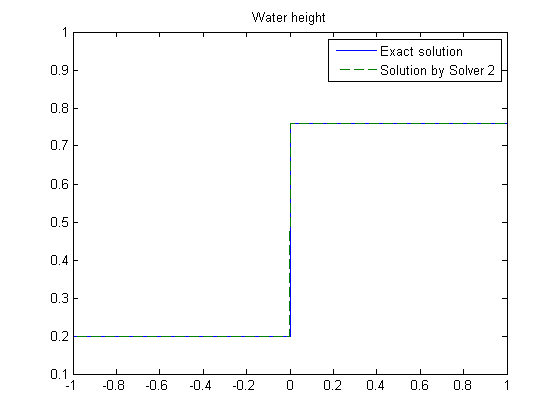}\\
  \caption{Test 6 (Resonant case - multiple solutions).
 Water height of the solution \rf{3.3} preferred Solver (A2).}\label{Resonance2h}
\end{figure}

\begin{figure}
  \includegraphics[width=10cm]{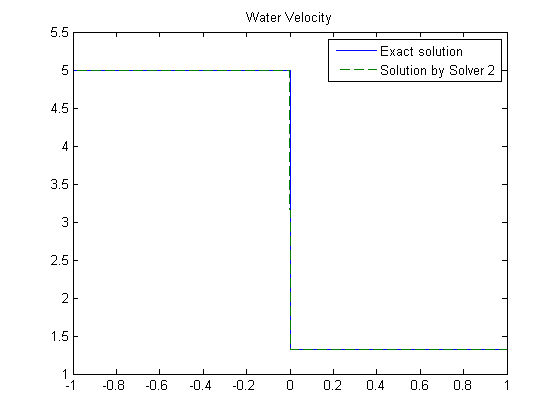}\\
  \caption{Test 6 (Resonant case - multiple solutions).
  Water velocity of the solution \rf{3.3} - preferred Solver (A2).}\label{Resonance2u}
\end{figure}

\begin{figure}
  \includegraphics[width=10cm]{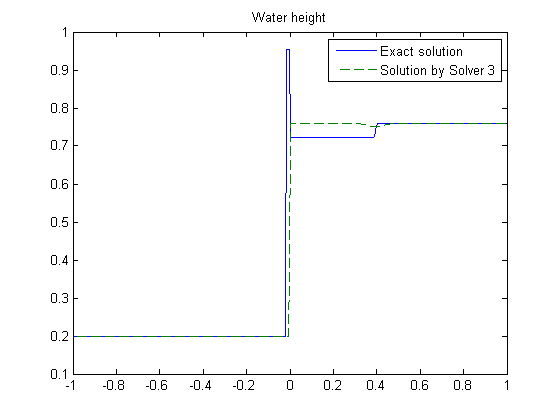}\\
  \caption{Test   (Resonant case - multiple solutions).
  Water height of the solution \rf{3.5} - preferred Solver (A3).}\label{Resonance3h}
\end{figure}

\begin{figure}
  \includegraphics[width=10cm]{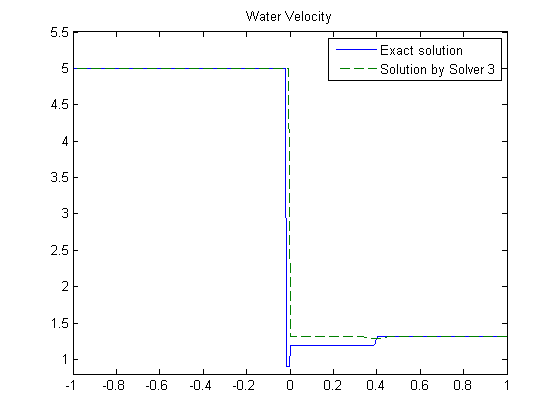}\\
  \caption{Test 6 (Resonant case - multiple solutions).
  Water height of the solution \rf{3.5} - preferred Solver (A3).}\label{Resonance3u}
\end{figure}

The Riemann problem \rf{1.1}--\rf{1.3} with initial data \rf{6.5} admits three distinct solutions: one solution of the form \rf{3.1} with
\begin{equation}
U_L^o=(0.21984063,4.5487497,1.2),\quad U_M=(0.7964266,1.4737915, 1.2),
\label{6.6}
\end{equation}
one solution of the form \rf{3.3} with
\begin{equation}
U_M=(0.75904946, 1.3174372, 1.2),
\label{6.7}
\end{equation}
which can be seen to be a stationary solution,
 and one solution of the form \rf{3.5} with
 \begin{equation}
U_M=(0.95328169, 0.89892673, 1),\quad U_M^o=(0.72279573, 1.1855776, 1.2).
\label{6.8}
\end{equation}
We could have three extreme choices of a Riemann solver for the Godunov method in this case. This can be seen easily by saying that we prefer a particular solver whenever it is available.
From Figures \ref{Resonance1h} and \ref{Resonance1u}  it follows 
that if Solver (A1) is preferred whenever it is available, then the  approximate solution approaches the exact solution. The same observation for Solver (A2);
 see Figures \ref{Resonance2h} and \ref{Resonance2u}. However, it is not the case for
 Solver (A3): approximate solutions do not converge to the exact Riemann solution by \rf{3.5};
 see Figures \ref{Resonance3h} and \ref{Resonance3u}.

 \subsection{Test 7}

Finally, we consider the Riemann problem with data
\begin{equation}
U_0(x)= \left\{\begin{array}{ll}
&U_L=(h_L,u_L,a_L)=(1, 2, 1.1)\in G_2,\quad x<0,\\
&U_R=(h_R,u_R,a_R)=(0.8,4,1)\in G_1,\quad x>0.\\
\end{array}\right.\label{6.9}
\end{equation}

\begin{figure}
  \includegraphics[width=10cm]{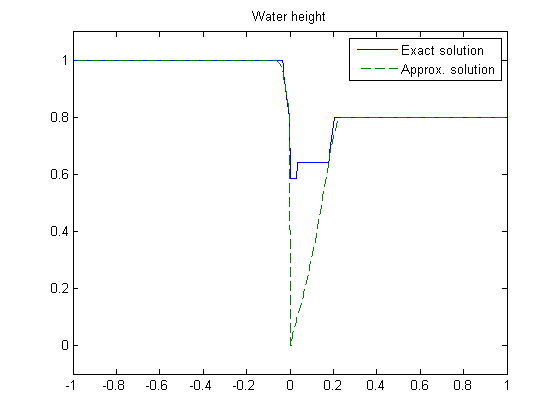}\\
  \caption{Test 7 (Resonant case).
 Water height of the solution  \rf{3.9}.
The Godunov scheme {\sl does
 not} approach the exact solution in this resonance case.}\label{B1-Resonance1h}
\end{figure}

\begin{figure}
  \includegraphics[width=10cm]{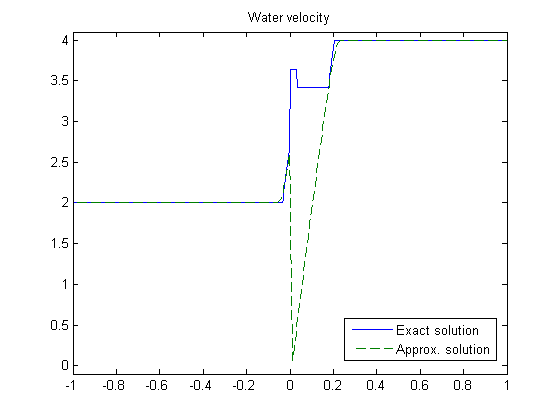}\\
  \caption{Test 7  (Resonant case.
 Water velocity for \rf{3.9}.
The Godunov scheme does not approach the exact solution in this resonance case.}\label{B1-Resonance1u}
\end{figure}
We have
$$
h_2^{\#}=  0.998204556070240  <h_1^o= 1.050890579855180,
  $$
  where $h_2^{\#}, h_1^o$ are defined as in  Theorem \rf{theo33}.
Thus, the problem \rf{1.1}, \rf{1.2}, \rf{6.1} admits a unique solution of the form \rf{3.9},  according to Theorem \rf{theo33}, where
\begin{equation}
\begin{array}{ll}
&U_1=(0.77374106, 2.7536634, 1.1), U_2=(0.58589019, 3.636556,1),\\
&U_3=(0.64142927, 3.4143821,1).
\end{array}
\label{6.10}
\end{equation}
Precisely, the solution begins with a $1$-rarefaction wave from $U_L$ to $U_1$, followed by a stationary contact from $U_1$ to $U_2$, then continued by a 1-shock wave from $U_2$ to $U_3$, and finally attains $U_R$ by a 2-rarefaction from $U_3$.
Figures \ref{B1-Resonance1h}-\ref{B1-Resonance1u}  show that  the Godunov scheme generates approximate solutions which {\sl do not} approach the exact solution in this resonance case (here shown at time $t=0.03$).

\section{Concluding remarks}

In this paper, we have provided for the first time a 
 complete characterization of all solutions to the Riemann problem associated with the shallow water equations.
First, we determined domains in the phase space in which
precisely one solution exists and domains in which several solutions (up to three) are available.
Second, we provided a computing strategy which allowed us to numerically determine 
all possible Riemann solutions.
Third, we defined a well-balanced and quasi-conservative Godunov scheme, which was
carefully tested in several regimes of interest.

The following main conclusions were established in this paper: 

\begin{itemize}

\item[$\bullet$] The proposed scheme captures exactly the equilibrium states.

\item[$\bullet$] It converges to the uniquely defined solution to the Riemann problem in strictly hyperbolic domains of the phase space, and this property validates our numerical strategy. We emphasize that 
the existing literature restricts attention to this strictly hyperbolic regime, only.

\item[$\bullet$]  Next, in order to further test our new algorithm, we
considered Riemann data belonging to both sides of the resonance curve;
both convergence to the selected solution and as well as convergence to another solution
were observed. This is not surprising since multiple solutions are available in such a regime.

\item[$\bullet$] Finally, the most challenging test was performed by taking
states precisely on the resonance curve, and we observed that
the scheme gave quite unsatisfactory results with no convergence whatsoever observed.
This latter behavior should not be interpreted as a drawback of the numerical method itself,
but rather indicates a limitation of the physical model which, in itself, does not properly describe the
fluid flow so that further physics is required in this regime. It is conceivable that a more satisfactory model 
could be obtained by analyzing the (small-scale) effects of higher-order terms and, for instance, introducing a suitable notion of kinetic relation, similarly to what is done for undercompressive shock waves \cite{LeFloch2002,LeFloch10}.

\end{itemize}

\end{document}